\UseRawInputEncoding
\documentclass[10pt]{amsart}
\usepackage{graphicx,amssymb,amsmath}
\newtheorem{thm}{Theorem}[section]
\newtheorem{cor}[thm]{Corollary}
\newtheorem{lem}[thm]{Lemma}
\newtheorem{prop}[thm]{Proposition}

\theoremstyle{definition}
\newtheorem{defn}[thm]{Definition}
\theoremstyle{remark}
\newtheorem{rem}[thm]{Remark}
\numberwithin{equation}{section}
\usepackage{amsmath,amssymb,amsthm}

\theoremstyle{plain}
\newtheorem{theorem}{Theorem}[section]

\theoremstyle{definition}

\def\Z   {{\bf Z}}
\def\Q   {{\bf Q}}
\def\R   {{\bf R}}
\def\C   {{\bf C}}
\def\P   {{\bf P}}
\def\F   {{\bf F}}

\def\Pic   {{\rm Pic}}

\def\HZ      #1#2{H^{#1}(#2,\Z)}

\def\HR      #1#2{H^{#1}(#2,\R)}

\def\calY      {{\mathcal Y}}
\def\calK      {{\mathcal K}}

\def\O           {{\mathcal O}}
\def\calE       {{\mathcal E}}
\def\calZ       {{\mathcal Z}}

\def\Kbar      {\overline{\mathcal K}}

\def\Kahler    {K\"ahler\ }
\def\CY        {Calabi--Yau\ }

\def\phimapX   {\phi : X \to \bar X }

\def\Mov       {{\it Mov}}
\def\Big        {{\it Big}}
\def\Movb     {{\overline \Mov}}
\def\Eff         {{\it Eff}}
\def\Effb       {{\overline \Eff}}
\def\Bigb      {{\overline \Big}}
\def\vol         {{\rm vol}}
\def\Cone      #1#2{{\rm Cone}  \langle #1 , #2 \rangle }
\def\Ttz         {{Type ${\rm III}_0$ }}

\def\Z   {{\bf Z}}
\def\Q   {{\bf Q}}
\def\R   {{\bf R}}
\def\C   {{\bf C}}
\def\P  {{\bf P}}
\def\F   {{\bf F}}

\newcommand{\beas}{\begin{eqnarray*}} 
\newcommand{\eeas}{\end{eqnarray*}} 

\begin{document}

\title{Boundedness questions for Calabi--Yau threefolds}



\author{P.M.H. Wilson} 
\address{Department of Pure Mathematics, University of Cambridge,
16 Wilberforce Road, Cambridge CB3 0WB, UK}
\email {pmhw@dpmms.cam.ac.uk}

\date{April 25, 2023}

\begin{abstract}
In this paper, we study boundedness questions for (simply connected) smooth Calabi--Yau threefolds.  The diffeomorphism class of such a threefold is known to be determined up to finitely many possibilities by the integral middle cohomology and two integral forms on the integral second cohomology, namely the cubic cup-product form and the linear form given by cup-product with the second chern class.  The motivating question for this paper is whether knowledge of these cubic and linear forms determines the threefold up to finitely many families, that is the moduli of such threefolds is bounded.  If this is true, then in particular the middle integral cohomology would be bounded by knowledge of these two forms.

Crucial to  this question is the study of rigid non-movable surfaces on the threefold, which are the irreducible surfaces  that deform with any small deformation of the 
complex structure of the threefold but for which no multiple moves on the threefold. 
 If for instance there are no such surfaces, then the answer to the motivating question is yes (Theorem 0.1).  In particular, for given cubic and linear forms on the second cohomology, there must exist such surfaces for large enough third Betti number (Corollary 0.2).  

The paper starts by proving general results on these rigid non-movable surfaces and boundedness of the family of threefolds.  The basic principle is that if the cohomology classes of these surfaces are also known, then boundedness should hold (Theorem 4.5).  The second half of the paper restricts to the case of Picard number 2, where it is shown that knowledge of the cubic and linear forms does indeed bound the family of Calabi--Yau threefolds (Theorem 0.3).  This appears to be the first non-trivial case where a general boundedness result for Calabi--Yau threefolds has been proved (without the assumption of a special structure).

\end{abstract}

\maketitle


\section*{Introduction}

\rm Let $X$ be a (simply connected) smooth complex \CY threefold.  A hard unsolved problem is whether such threefolds form a bounded  (or even 
just birationally bounded) family.   
This would in turn imply that the Euler characteristic for \CY threefolds is bounded. 
For example, by results of Gross \cite{GrEll}, elliptically fibred \CY threefolds do form a birationally bounded family; no such result 
is known for fibre spaces over $\P^1$ with generic fibre a K3 or abelian surface.
One can split the general problem into two parts; 
whether or not \CY threefolds fall into a finite 
number of topological types, and whether or not  \CY threefolds of a given topological type form a bounded family.  The first of these 
problems seems intractable whilst the second is hard but maybe tractable.
It is this latter problem which largely motivates much of the theory developed in this paper.

The diffeomorphism class of $X$ is determined 
up to finitely many possibilities by knowledge of certain topological invariants, namely
the cup-product cubic form on $H^2 (X, \Z)$ given by $D \mapsto D^3$, the linear form 
on $H^2 (X, \Z)$ given by $D \mapsto D\cdot c_2(X)$ and the middle cohomology $H^3 (X, \Z)$ \cite{Sull}, 
and if furthermore $H_2 (X,\Z )$ is torsion free, this
information determines the diffeomorphism class precisely \cite{Wall}.  

If the linear form $c_2$ were trivial, it would follow from Yau's proof of the Calabi Conjecture \cite{Yau} that the Ricci flat metric on $X$ would in fact be flat, implying that $X$ is an \' etale quotient of an abelian threefold, 
contradicting our assumption that $X$ is  simply connected.

In this paper, we address the 
question  
as to whether $X$ is determined up to finitely many families by knowledge of the cubic and linear forms on $H^2 (X, \Z)$; if this is true then of course the diffeomorphism type, and in particular $H^3 (X, \Z)$, will have only finitely many possibilities.  This would contrast with the classical results of C.T.C. Wall on the diffeomorphism types of 6-manifolds \cite{Wall}, which imply that for any given allowable data of cubic and linear forms on 
$H^2 (X, \Z)$, the value of the third Betti number is unbounded, since one can always take connected sums with an 
arbitrary number of copies of $S^3 \times S^3$.  Even for (non-K\"ahler)
complex Calabi--Yau threefolds which admit balanced metrics (i.e. $d(\omega ^2 ) =0$), a similar flexibility occurs and there 
exist for instance examples with $b_2 =0$ but with $b_3$ arbitrarily large \cite{FLY}.

For the above boundedness question,  a major role will be played by irreducible surfaces $E$ on the \CY threefold $X$ 
that deform with any small deformation of the 
complex structure of the threefold but for which no multiple moves on the threefold. 
 We use the terminology \it rigid non-movable surfaces \rm for these 
(Definition 2.1).  We will note in Section 2 that there is a birational description of such surfaces and that they all 
 contain at least a one dimensional family of rational curves.
The main results proved in this paper on the stated boundedness question are the following 
three general results

\begin{thm}  For \CY threefolds $X$ containing no rigid non-movable surfaces, knowledge of the cubic cup-product form and the linear form $c_2$ on 
$H^2 (X, \Z)$ ensures that $X$ lies in a bounded family.\end{thm}

\begin{cor} For \CY threefolds $X$ with given cubic and linear forms on $H^2 (X, \Z)$,  there must exist rigid non-movable surfaces on $X$ when  $b_3(X) \gg 0$ .\end{cor}

In fact, as explained in Remark 4.7, a stronger result than Corollary 0.2 holds.  The corollary might be compared with the main result from \cite{HBW}, where the existence of rational curves is shown 
when $b_2(X) > 13$.

In the case of Picard number $\rho =1$, there will be no rigid non-movable surfaces.  Moreover the generator $L$ of $\Pic (X)$ with $L^3 >0$ will be ample, and so 
by a simpler version of the Hilbert scheme argument from Proposition 1.1, knowledge of the cubic and linear forms on $H^2 (X, \Z)$ does ensure that $X$ lies in a bounded family 
 --- in fact 
it can be seen that we only need the 
cubic form for this.  In this paper we can therefore assume $\rho >1$ throughout.

For higher $\rho$, the interplay between the possible rigid non-movable surface classes is rather delicate, and for $\rho > 2$ there are only partial results
on boundedness.
For $\rho =2$, where closed convex cones 
in $\HR 2 X$ are determined by their edge rays, we can do better.

 \begin{thm}   For \CY threefolds with Picard number $\rho (X) =2$,  
 knowledge of the cubic cup-product form and the  linear form $c_2$ on 
$H^2 (X, \Z)$ ensures that $X$ lies in a bounded family.\end{thm}

As far as the author is aware, this 
is the first non-trivial case where a general boundedness result for \CY threefolds has been proved (without the assumption of a special structure).
The proof of this  theorem throws up a number of interesting questions about properties which hold for $\rho =2$ and that are less clear for higher Picard number.

\section{Preliminaries}

In this section, we recall and elaborate on various results from the literature 
for \CY threefolds, largely revolving around the various cones of divisor classes 
contained in $\HR 2 X$ and the change in the above cubic and linear forms on 
$\HR 2 X$ under flops, and we explain the underlying philosophy behind the proofs of the 
main results.  We note in passing that, since $X$ is assumed simply connected, the abelian group 
$\Pic (X) \cong \HZ 2 X$ is torsion-free, and so the integral elements of $\HR 2 X$ do correspond 
precisely to the divisor classes.  As explained in the Introduction, we shall assume that the Picard number $\rho >1$.

If $X$ is a smooth \CY threefold with \Kahler cone $\calK$, then the
nef cone $\Kbar \subset \HR 2 X$ is locally rational polyhedral away from the cubic cone 
$$W^* = \{ D\in \HR 2 X \ : \ D^3 =0\};$$ moreover, the codimension one faces of
$\Kbar$ (not contained in $W^*$) correspond to primitive birational contractions 
$\phimapX$ of one of three different types \cite{WilKC}.

In the numbering of \cite{WilKC}, Type I contractions are those where only a finite
number of curves (in fact $\P ^1$s) are contracted.  The singular threefold
$\bar X$ then has a finite number of cDV (compound Du Val) singularities.  Whenever one has such a primitive 
small contraction on $X$, there is a flop of $X$ to a different birational
model $X'$, also admitting a birational contraction to $\bar X$; moreover,
identifying $\HR 2 {X'}$ with $\HR 2 {X}$, the 
nef cone of $X'$ intersects the nef cone of $X$ along the codimension one face
which defines the contraction to $\bar X$ \cite{Kaw, Koll}.  It is well known (e.g. \cite{Koll}, Theorem 5.2.3)
that $X'$ is smooth, projective and has the same Hodge numbers as $X$, but that
the finer invariants, such as the cubic form on $\HZ 2 X$ given by cup-product,
and the linear form on $\HZ 2 X$ given by cup-product with 
the second chern class  $c_2 (X) = - p_1 (X)/2$,
will in general change.  

We shall call  an integral divisor class on a smooth complex projective variety  \it mobile \rm if the corresponding linear system is non-empty and  has no fixed 
components.  Taking the closure of the cone generated by mobile classes yields a closed 
cone  $\Movb (X)  \subset \HR 2 X$, which we shall refer to as the \it movable \rm cone, although 
this term is sometimes elsewhere applied to its interior $\Mov (X)$.  We shall call a class in the closed  cone 
\it movable\rm , and a class in the open cone \it strictly movable\rm .  

For $X$ a \CY threefold, a result of Kawamata (\cite{Kaw},  page 120)
decomposes the open cone $\Mov (X)$ into chambers, each corresponding to the proper transform of 
the \Kahler cone of some birational minimal model of $X$ under a sequence of flops, 
a codimension one  wall between adjacent chambers corresponding to the flop between the corresponding minimal models, 
thus defining a Type I contraction on both of them.

The other cone of divisor classes which will be of interest to us is the closure of the cone generated by 
effective classes. the \it pseudo-effective \rm cone  $\Effb (X)$, and its interior the \it big \rm cone  $\Big (X)$; hence 
$\Bigb (X) = \Effb (X)$.
Thus $ \Movb (X) \cap \Big (X)$ consists of big movable classes; a strictly movable class is big and movable, 
although the converse does not in general hold.   

For \CY threefolds, any rational class in $\Movb (X) \cap \Big (X)$ 
corresponds 
(under a sequence of flops) to a  big nef class on some birationally equivalent minimal model (see proof of Theorem 5.3 in \cite{Kaw}).  In fact, 
any rational class in $\Big (X)$ has an integral multiple which can be written as a mobile divisor plus an effective divisor,
where the mobile divisor defines a birational map.  The advantage of the convex cones 
just defined is that, unlike the \Kahler cone,  they are 
birational invariants.

Another crucial tool that we shall constantly use is the formula for how the cubic form and the linear form 
$c_2$ transform under flops from one 
minimal model to another.  This is explained in \cite{WilGW}, where it is observed that by locally deforming the complex structure in 
a neighbourhood of the exceptional locus of the Type I contraction, we may reduce down to the case of a disjoint set of $(-1, -1)$-curves, which may not necessarily
be achievable by a global complex deformation of $X$ (although will be achievable by an almost complex deformation).  

The idea behind this reduction may be found on page 679 of \cite{FrSR}: the contraction of a connected component $Z$ of the exceptional locus is a cDV (compound Du Val) singularity, for which we 
can take a Stein neighbourhood $\bar U$, with corresponding open set $U \supset Z $ in $X$.  A generic 
section of $\bar U$ through the singularity is then a rational double point $Y_0$, 
and we have a holomorphic map $g : \Delta \to \hbox{Def}\, Y_0$ to the versal deformation space of the singularity.  Pulling back 
the flat family  to $U$, we have a partial resolution $X_0$ 
of $Y_0$ and $g$ lifts to a holomorphic map $f : \Delta \to \hbox{Def}\,  X_0$, the versal deformation space of the partial resolution, where there is a natural morphism 
$\pi : \hbox{Def}\, X_0 \to \hbox{\rm Def}\, Y_0$.  

If $\tilde Y_0$ is the minimal resolution of $Y_0$, then $\hbox{Def}\, Y_0$ is the quotient of $\hbox{Def}\,  \tilde Y_0$ by the 
Galois action of the Weyl Group $W$ of the singularity, and $\hbox{Def}\, X_0$ is the quotient of $\hbox{Def}\,  \tilde Y_0$ by a subgroup  
$W_0$, the Weyl group  associated to 
the $(-2)$-curves on $\tilde Y_0$ contracted under the map to $X_0$ (see Section 8 of \cite{Pink}).  If the rational double point singularity has rank $n$, then all these spaces are just neighbourhoods of the origin in 
$\C ^n$, and 
on each deformation space we have a (reduced) divisor given by the inverse image of the discriminant locus on 
$\hbox{Def}\, Y_0$, which on $\hbox{Def}\, X_0$ we denote by $D$, and on $\hbox{Def}\, \tilde Y_0$ is a collection of hyperplanes (through the origin) 
corresponding to pairs $\{ \pm {\bf r}\} $ for ${\bf r}$ a 
root of the system (so that $D$ is the image of these hyperplanes under the quotient map).  Deforming $f$ so that it is transverse to $D$ gives the local 
deformation of $U$ under which $Z$ splits up into a finite number $\delta$ of $(-1, -1)$-curves (cf. \cite{WilSD}, Section 1). 
Given such a deformation of $f$, for each irreducible component $D_i$ of $D$ we have that the deformation is transverse to $D_i$ at points not on any other component of $D$.  So each irreducible component of $D$ contributes at least one to $\delta$.
In particular, we can explicitly calculate the number $r$ of irreducible 
components of $D$ from the action of $W_0$ on the roots, and then we may observe $\delta \ge r$. This inequality plays a crucial role in some of the calculations below.

The simplest case here is when $Y_0$ is an $A_1$-singularity, where $r=1$ and $\delta$ is the ramification index of $f$, and we have the 
situation studied in \cite{Reid}.  Conversely, by considering the action of $W_0 \subset W $ on the roots, if $Y_0$ is not an $A_1$-singularity we see that $r>1$.  In particular we deduce that $\delta =1$ if and only if $Z$ itself is a $(-1,-1)$-curve.

The next case is when $Y_0$ is an $A_2$-singularity; here 
the results of \cite{KatzM, KawHyp} imply that $Z$ cannot be irreducible since otherwise the generic section of $\bar U$ would be an $A_1$-singularity.
Hence we have $X_0 = \tilde Y_0$,  in which case $Z$ has two irreducible components and $r=3$, 
corresponding to three lines through the origin.  
Therefore $\delta \ge 3$ 
with equality only if $f$ is locally an isomorphism and transverse to each of the three lines.  

It may be checked 
from the explicit description of roots and Weyl groups found in say Section 2 of \cite{KatzM} that if $Z$ has precisely two components and we are not in the specific $A_2$ case described above (so either $Y_0$ is an $A_2$-singularity but $f$ is ramified or not generic, or  $Y_0$ is worse than an $A_2$-singularity and so $X_0$ is only a partial resolution of $Y_0$), 
then $\delta \ge 4$ (in the latter case of $Y_0$ worse than an $A_2$ singularity we have  $r\ge 4$).  If for instance $Y_0$ is an $A_n$-singularity ($n>2$), and we are contracting all 
$(-2)$-curves on $\tilde Y_0$ except for those  corresponding to two adjacent nodes in the Dynkin diagram, then in fact $r=4$.

Knowledge of these 
\it virtual \rm $(-1, -1)$-curves is enough to determine the transformation law for the cubic and linear forms.  Let $\phi : X \to \bar X$ be the Type I contraction, 
and $\eta \in H_2 (X , \Z)/{\hbox{\rm Tors}}$  the primitive class contracted by $\phi$, i.e. $\phi _* \eta =0$.  Let $n_d$ denote the total number of $(-1, -1)$-curves in the above 
local deformations (corresponding to the various connected components of the exceptional locus) which have image the class $d\eta$ in $H_2 (X, \Z)/{\hbox{\rm Tors}}$.  
It is shown in \cite{WilGW} (building on the theory from \cite{WilSD} and \cite{Fr}) that the cubic and linear forms transform according to the following formulae:
$$
 D_2^3=D_1^3-(D_1\cdot\eta)^3\sum_{d>0}n_d d^3, \ \  
 D_2 \cdot c_2(X')=D_1\cdot c_2(X)+2(D_1\cdot\eta)\sum_{d>0}n_d d \  \ 
  \eqno (1)$$
 where $X'$ denotes the flopped threefold and $D_2$ is the divisor on $X'$ corresponding to a divisor $D_1$ on $X$.  This then leads to a pivotal observation.
 
 \begin{prop} For \CY threefolds $X$, knowledge of $D^3$ and $D\cdot c_2$ for some 
  big movable  class 
$D \in H^2 (X, \Z)$ determines $X$ up to a bounded family.
 
 \begin{proof}
 
 Suppose we have a big movable integral divisor class $D$ on $X$, by 
 Section 5 of \cite{Kaw} or \cite{KollarFlops} 
 Corollary 6.3 therefore corresponding to 
 a big nef class $D'$ on some birationally equivalent (smooth) minimal model $X'$; note that $c_2 (X')\cdot D'$ is non-negative \cite{Yau, Miy}.  
 The birational transformation from $X$ to $X'$ is obtained 
 by successively making directed flops in curves on which $D$, and subsequently the transforms of $D$, are negative, also called a sequence of $D$-flops  \cite{Kaw, KollarFlops, KM2}.  Thus successively applying the formula (1)
 for the transformation of the linear form,  
 we obtain $c_2(X)\cdot D$ in terms of $c_2(X')\cdot D'$ with 
 $c_2(X)\cdot D \ge c_2(X') \cdot D' \ge 0$, and so knowledge of $c_2 (X)\cdot D$ bounds all the additional terms on the right-hand sides of both the equations (1).  If furthermore we 
 know $D^3$, then the transformation law for the cubic form will then bound $(D')^3$.  Thus knowing a big movable integral class $D$ on $X$ will yield a big nef integral 
 class $D'$ on a birationally equivalent minimal model $X'$, with bounds on both $c_2 (X')\cdot D'$ and  $(D')^3$.  When we factor the above birational map into  
 a sequence of $D$-flops $X = X_m -\to X_{m-1}  \ldots -\to X_1 -\to X_0 = X'$, we denote by $D_i$ the big mobile divisor on $X_i$ corresponding to $D$ on $X$ (thus $D_0 =D'$), and note that $m$ is bounded by $ c_2(X)\cdot D /2$.

 If $D'$ is ample, we know that $10D'$ is very ample \cite{OP}; if $D'$ 
 big but not ample, 
 then it determines a birational contraction $\phi : X' \to Y$, with $Y$ a singular \CY threefold with canonical singularities, and $D' = \phi ^* L$ for some ample Cartier divisor $L$ on $Y$.  The results of \cite{OP} then imply that $14 L$ is very ample on $Y$.  Applying the theory of Hilbert schemes, we deduce in both cases that $X'$ lies in a bounded family (in the 
 second case, coming from strata of a Hilbert scheme over which the corresponding threefolds are equisingular \CY with canonical singularities, under crepant resolution of 
 the singularities).  Note here that a \CY threefold with canonical singularities only has finitely many possible crepant resolutions by 
 \cite{KollarFlops}  Corollary 5.6 or  \cite{KM2} Theorem 6.42.
 Knowledge of the values of the cubic and linear forms on 
 a big movable class on $X$ therefore restricts $X$ to lie in a \it birationally \rm  bounded family.
 
 We say that $X'$ (and hence each possible $X_i$) is \it general \rm in moduli if the only non-movable irreducible surfaces contained in it are rigid; 
 the locus of such general moduli points 
consists of  the complement of countably many proper subvarieties.  Since the corresponding moduli space of (polarised) \CY threefolds is smooth, it may be shown that between any two general points in moduli, one can find a piecewise continuously differentiable simple real curve 
of points general in moduli (the subvarieties are of real codimension at least two and we use the Baire category theorem locally on coordinate neighbourhoods).  Over such a curve it follows from \cite{WilKC} that the \Kahler cone is constant (with any flopping face remaining so in the family), and that 
moreover we may simultaneously flop 
in such a family by \cite{KM}, Theorem 11.10.  This reduces us to showing that starting from a known  $X'$ general in moduli, there are only finitely many sequences 
of possible flops satisfying the above numerical bounds, and hence only finitely many  possible  families for $X$.

 For a given $0\le i\le m-1$, 
we suppose inductively that the composition of flops $X_0 -\!\to X_i$ has been determined up to finitely many possibilities, and we show that for each possibility there can be be only finitely 
many possibilities for the next flop $X_i -\!\to X_{i+1}$.  We choose a very ample divisor $H_i$ on $X_i$, with corresponding (known) mobile and big divisor $H_0$ on $X_0$ and 
 divisor $H_{i+1}$ on the yet to be determined minimal model $X_{i+1}$.  Since $X_0$ is general in moduli, we can find a positive integer $n$ such that $nD' -H_0 = \Delta 
+\sum_{i=1}^r s_j E'_j$, where $\Delta$ is movable on $X_0$, the $E'_j$ are (known) rigid non-movable surfaces in the sense of Definition 2.1 
and the $s_j$ are (known) non-negative integers, all depending on our choice of $H_i$.  If $E_1, \ldots ,E_r$ denote the corresponding rigid non-movable surfaces on $X_{i+1}$, we deduce that $nD_{i+1} - H_{i+1} - \sum_{i=1}^r s_j E_j$ is movable on $X_{i+1}$, and this divisor  will then,   as noted in Remark 1.2, have non-negative intersection with $c_2 (X_{i+1})$.  Thus $c_2(X_{i+1})\cdot H_{i+1} + \sum s_j c_2(X_{i+1})\cdot E_j \le nc_2(X_{i+1})\cdot D_{i+1}$.  
Now the formula (1) implies that $c_2(X_k)\cdot D_k \le \ c_2(X)\cdot D$ for all $0\le k\le m$.  Moreover, since the $E_j$ are rigid non-movable surfaces, it follows 
from the arguments of Section 2 below, and in particular equation ($1'$) and inequalities (3), that $c_2 (X_{i+1})\cdot E_j \ge -6$ for all $j$.  Thus we deduce that 
$c_2 (X_{i+1})\cdot H_{i+1} \le n c_2(X)\cdot D + 6\sum s_j$, an upper bound dependent only on known data.  We now use the formula (1) to relate $c_2 (X_{i+1})\cdot H_{i+1}$ and $c_2(X_i)\cdot H_i$;  the flop $X_{i +1} -\!\to X_i$ is an $H_{i+1}$-flop and the sum of the additional virtual terms appearing on the right-hand side of the formula is bounded, and this in turn  implies a bound (depending on known data) on 
$H_i \cdot C$ for any irreducible curve $C$ in the flopping locus of $X_i -\!\to X_{i+1}$.  This restricts the rigid curve $C\cong \P^1$ to being one of a finite number of curves on $X_i$ by considering the corresponding Hilbert scheme of such curves of bounded degree with respect to $H_i$, and this in turn implies that there are only finitely many possibilities for the flop $X_i -\!\to X_{i+1}$.  The proof of boundedness is therefore complete.
\end{proof}\end{prop}

So given a big movable integral class $D$, Proposition 1.1 implies that there are only finitely many possibilities for the Euler number $e(X)$.  This can however be made 
explicit in the case when $D$ corresponds to an ample divisor $H$ on a 
birationally equivalent (smooth) minimal model, i.e. when $D$ lies in the interior of one of the chambers 
in the decomposition of  $\Mov (X)$.
There is an explicit bound on $H^3$, for $H$ ample corresponding to $D$ on the flopped model, a function of $D^3$ and $D\cdot c_2$.  However $| e(X) |$ is then bounded 
by an explicit multiple of $H^3$ by results in \cite{KW}, and hence by a function of $D^3$ and $D\cdot c_2$. 
  
 \begin{rem}  The argument in the first part of the proof of Proposition 1.1 shows that $c_2 (X)\cdot D \ge 0$ for any strictly movable rational class $D$ on a \CY threefold $X$, and hence that the linear form $c_2 (X)$ is non-negative on $\Movb (X)$.  This remark plays a crucial role in the proofs in the last two sections of this paper. \end{rem}
  
Given knowledge of the cubic and linear forms on $H^2 (X, \Z)$, the aim 
throughout this paper will be to find a finite set of integral classes, determined by given cohomological data but maybe involving 
a finite number of  arbitrary choices for integral classes in certain  cones determined by the data, where for 
a given \CY threefold $X$ general in moduli with this data, one of these classes 
must represent a big movable divisor.  Thus Proposition 1.1 will then  prove boundedness for all 
\CY threefolds with the given cohomological data.
Crucial to this problem will be the rigid non-movable surfaces contained in $X$.  

 \section{Rigid non-movable surfaces}

Given an irreducible surface $E$ on a \CY threefold $X$, we may perform a sequence of directed flops on $X$ so that either $E$ corresponds to a nef divisor on a birationally equivalent 
minimal model $X'$, or the corresponding surface $E'$ on $X'$ may be contracted 
(\cite{Kaw} Theorem 7.1, \cite{KollarFlops} Corollary 6.3).  
In the former case, we know that some multiple of the (irreducible) nef divisor moves and defines a  morphism on $X'$ \cite{Og}.  In
the latter case,  the birational contraction is either of Type II, in which case $E'$ is a generalised del Pezzo surface, 
or it is of Type III and $E'$ is a conic bundle over a smooth curve $C$ of genus $g$; see \cite{WilSD} for further details.  The cases of Type II contractions and Type III contractions with $g=0$, which from now on will be denoted Type ${\rm III}_0$, will be of particular interest in this paper.

Given the previously noted result which enables us to flop in families
(\cite{KM}, Theorem 11.10), the surface $E$ deforms 
with small deformations of the complex structure on $X$ if and only if $E'$ deforms under small deformations of the complex structure on $X'$.  The exceptional surface of a Type II contraction on $X'$ always deforms as a surface (cf. proof of Proposition 3.2(ii) in \cite{WilKC}), 
and the same is therefore true of the corresponding surface $E$ on $X$.  If $E'$ is the exceptional surface of a Type III contraction, it will not deform as a surface if $g>0$ but will deform if $g=0 $
(see \cite{WilSD}, \S 4), and this then determines whether $E$ deforms 
or not under small deformations 
of the complex structure on $X$.  In the cases when $g>0$, the Type III contraction becomes a Type I contraction on a generic small deformation, or in the special cases studied in \cite{WilKC} ceases to be a contraction.

For proving boundedness of families, we can assume the threefold is general in moduli, 
i.e. in the complement of a countable union of subvarieties in its moduli space, 
and so we 
will be able to ignore those surfaces which do not deform under generic small deformations.  We are therefore led to a basic definition.

\begin{defn} If $X$ denotes a Calabi--Yau threefold, a \it rigid non-movable surface \rm on $X$ is an irreducible surface $E$, whose divisor class in $H^2 (X,\Z )$ represents a surface on any small deformation of the complex structure on $X$, but for which no integer multiple is mobile. 
In particular, if $X$ is general in moduli, any irreducible surface on $X$ is either movable or rigid non-movable. 
\end{defn}

Thus from the above discussion, the rigid non-movable surfaces on $X$ are precisely those which correspond on some minimal model
to the exceptional surface of a contraction of Type II or Type ${\rm III}_0$.

We now recall known results concerning such exceptional surfaces.  First we consider the case when $E$ is the exceptional surface of a Type II contraction, where we quote from the discussion in 
\cite{WilSD}, \S 2.  Here $E$ is an irreducible \it generalised del Pezzo surface, \rm as classified in Theorem 1.1 of \cite{ReiddP}. The normal 
irreducible generalised 
del Pezzo surfaces are either elliptic cones or del Pezzo surfaces with rational double point singularities (cf. page 620 of \cite{WilSD}). 
By considering the local embedding dimension at the vertex, the only elliptic cone which can occur in our case is a  cone in $\P^3$ on a smooth plane cubic, where $E^3 =3$. For $E$  a del Pezzo surface with 
rational double point singularities, we have 
$1 \le E^3 \le 9$, where $E\cong \P^2$ if $E^3 =9$, and 
$E \cong \F_1$, $\P^1 \times \P^1$ or a quadric cone if $E^3 =8$, where $\F_1$ denotes the Hirzebruch surface
(namely the blow-up of $\P^2$ in a point).  However 
the only case 
of a Type II contraction with $E^3 >3$ where $E$ might be non-normal is when $E^3 =7$ and 
$E$ is one of the surfaces $\bar\F_{3,2}$ or $\bar\F_{5,1}$ described in the proof of 
Lemma 2.3 in \cite{WilSD} --- see also the proof of Theorem 5.2 in \cite{GrDef}, noting also Remark 5.3 there.  These specific non-normal surfaces both contain a curve of double points isomorphic to $\P^1$, the generic section of which is a simple node. Finally we recall that all the surfaces $E$ occurring here are smoothable to  smooth del Pezzo 
surfaces --- this is clear when $E^3 \le 3$, 
the cases where the contraction is a hypersurface singularity (cf. pages 620-1 of \cite{WilSD}), 
and follows for all other cases by Lemma 5.6 of \cite{GrDef}.  In particular we have $\chi 
(E, \O_E) =1$ (cf. also the top of page 567 of \cite{WilKC}).  Hence, using the 
Riemann--Roch theorem, we get  $$1 = \chi (\O_E ) = 
\chi (\O_X) - \chi (\O_X(-E)) = \chi (\O _X (E)) = \frac{1}{6} E^3 + \frac {1}{12} c_2(X)\cdot E,$$ 
and thus 
$$2 E^3 +  c_2(X)\cdot E =12. \quad
  \eqno (2)$$ 
The moral therefore is that we have fairly precise information on the exceptional surfaces of Type II contractions.

We now consider the case when $E$ is the exceptional locus of a Type III contraction 
of $X$; $E$  is therefore a conic bundle over a smooth curve $C$ of genus $g$.  If $g >0$, 
then curves from 
only finitely many fibres will deform under a generic small deformation of the complex structure on $X$; if however $g=0$, then the whole divisor will always deform (see \cite{WilSD}, \S 4).  In other words, 
for $g=0$ the Type III contraction remains Type III under deformations, whilst for $g>0$ it either becomes a Type I contraction under generic deformations, or in the special case of 
$E$ an elliptic ruled or quasi-ruled surface 
studied in \cite{WilKC}, ceases to be a contraction.  When $g=0$, the exceptional surface $E$ might be non-normal, therefore having generic fibre over $C$ a line pair.  However under a generic deformation of $X$, the surface $E$ will deform to a normal surface over $\P^1$ (having irreducible generic fibre), except for the case when $E^3 =7$ and the deformation has generic fibre a line pair and precisely two double fibres over the base curve (\cite{GrPrim} Proposition 1.2 and Theorem 1.4,  and \cite{WilSD} Proposition 3.2); 
using the latter Proposition,  it was pointed out at the end of the proof of Proposition 4.2 of \cite{WilSD} that  $E$ is in fact the non-normal del Pezzo surface $\bar\F_{3,2}$ we encountered in  
the Type II case.  \it Assuming \rm $X$ is general in moduli, we therefore have a characterization of possible exceptional surfaces for \Ttz contractions, where a full description of the ones with 
irreducible generic fibre may be found in Lemma 3.1 of \cite{WilGW}.  In particular we deduce in general that $E^3 \le 8$ (although in this case it may be negative).  If the conic bundle is denoted by 
$\pi : E \to C$, then $R^1\pi _* \O _E = 0$ and $\pi _* \O _E = \O _C$ and so $\chi (E, \O _E ) = 
\chi (C, \O _C ) = 1$ as for the Type II case; in particular the above argument via the Riemann--Roch theorem 
shows that equation (2) holds true also for the \Ttz case.
Since $E^3 \le 9$ for the Type II case and $E^3 \le 8$ for the Type ${\rm III}_0$ case, we note therefore using equation (2) that 
$$c_2 (X)\cdot E \ge -6,  -4 \eqno(3)$$ 
respectively for the Type II case and  the \Ttz  case.  We noted above that for a precise classification 
of these exceptional surfaces, we should assume that $X$ is general in moduli, but  
equations (2) and (3) remain 
true for Type II or \Ttz exceptional surfaces without that assumption.

Suppose now that $E$ is any rigid non-movable surface on a \CY threefold $X$; it corresponds to an exceptional surface $E'$ of either a Type II or a \Ttz contraction 
on some birationally equivalent minimal model $X'$ (obtained from $X$ by a succession of directed flops).

In this situation, we observe from the formulae (1) in Section 1 that with the notation as there
$$
 E^3=(E')^3-(E'\cdot\eta)^3\sum_{d>0}n_d d^3, \ \ 
 E \cdot c_2(X)=E'\cdot c_2(X')+2(E'\cdot\eta)\sum_{d>0}n_d d \ \  \eqno (1')
 $$
or an iteration of such formulae if there is more than one directed flop involved.  At each stage, the flop is directed by the surface corresponding to $E$, so that 
the correction terms on the righthand side of the first equation in ($1'$) are always negative and the correction terms on the 
righthand side of the second equation  are always positive.  
We shall in what follows refer to the elementary flop in a $(-1,-1)$-curve as an \it Atiyah \rm flop.  From an observation in Section 1, these occur precisely when there is a unique non-zero $n_d$ and this 
has value one.
From the above  observations 
we deduce the following.

\begin{prop}   Let $E$ denote a rigid non-movable surface on a \CY threefold $X$.  
\smallskip

(i) If $c_2 (X) \cdot E $ is bounded above, then there are only finitely 
many possibilities for $E^3$.  
If $c_2 \cdot E < 0$, then all these possibilities have 
$E^3 >0$.  
\smallskip

(ii) If $E^3$ is bounded below, then 
there are only finitely many possibilities for $c_2(X)\cdot E$.  
\smallskip

(iii) Suppose now $\rho = 2$ and the cubic and linear forms on $H^2 (X, \Z)$ are known.  If 
the values  $E^3$ and $c_2 (X)\cdot E$ are specified, then there are only finitely many possible 
associated cohomology classes in $H^2 (X, \Z)$.  

\begin{proof}
(i) With $E' \subset X'$ defined as above, since $c_2 (X')\cdot E'$ is bounded below from equation (3), knowledge of an upper bound for $c_2(X)\cdot E$ restricts the 
correction terms on the right-hand sides of ($1'$) to a finite number of possibilities.  But $(E')^3$ is bounded below as $c_2 (X')\cdot E' \le c_2(X)\cdot E$ is bounded above and equation (2) 
holds for $E'$.
Thus the previous observation on the correction terms in ($1'$) implies that there are only finitely many possibilities for $E^3$.  

If $c_2 (X) \cdot E < 0$, then from equation (3) 
(as it is an even number)  it is $-6$, $-4$ or $-2$.  In the case of $c_2 (X) \cdot E = -6$, we know 
that $E$ itself is the exceptional surface of a Type II contraction and $E^3 =9$, namely $E \cong \P^2$. 
 If $c_2 (X) \cdot E 
= -4$, then either $E$ itself is the exceptional surface of a Type II or \Ttz contraction and $E^3 =8$, or we are making a single Atiyah flop going from $X'$ to $X$, where $c_2 (X') \cdot E' = -6$ and $E'\cdot \eta =1$, 
with $n_1 =1$ and $n_d = 0 $ for $d>1$.  In this case $(E')^3 =9$ and $E^3 =8$.

If $c_2 (X) \cdot E = -2$, there are three possibilities.  Firstly $E$ itself could be the exceptional surface of a Type II or \Ttz contraction with $E^3 =7$.  Secondly we might have $c_2 (X') \cdot E' 
= -4$ and a single Atiyah flop going from $X'$ to $X$ 
with $E'\cdot \eta =1$, 
$n_1 =1$ and  $n_d = 0$ for $d>1$.  In this case $(E')^3 =8$ and again $E^3 =7$.  Finally we might have $c_2 (X') \cdot E' 
= -6$, i.e. $E' \cong \P^2$; here there could be two flops from $X'$ to $X$, but both of these would be Atiyah flops in curves with intersection 
multiplicity one with the surface corresponding to $E$, 
and hence again $(E')^3 = 9$ and $E^3 =7$.  The remaining possibilities with $c_2 (X) \cdot E = -2$ and 
$c_2(X')\cdot E' = -6$ have \it either 
\rm $E' \cdot \eta =2$, $n_1 =1$ and $n_d =0$ for $d>1$, \it or \rm $E' \cdot \eta =1$, $n_2 =1$ and $n_d =0$ for 
$d \ne 2$; in both these cases we get $E^3 = 1$.

(ii) Since $(E')^3$ is bounded above (by 9 for the Type II case, and 8 for the Type ${\rm III}_0$ 
case), a similar argument to that for the first part of (i), using the equations ($1'$), (2) and (3),  
shows  that there are only finitely many possibilities for $c_2(X)\cdot E$.

(iii)   This part follows immediately, 
unless $E^3 =0$ and $c_2(X)\cdot E =0$, from the fact that the linear form $c_2$ is non-trivial and 
an affine cubic in one variable has 
at most three roots.  In the remaining case we observe that the condition $c_2 \cdot E=0$  defines a rational line 
through the origin.  
We then note from equations ($1'$) 
that there must be some flopping curve $C$ on $X$, in fact a $(-1, -1)$-curve, with $E\cdot C = -1$ or $-2$, and therefore only finitely many integral classes within the line can occur.
\end{proof} \end{prop}

We comment that if we wish to extend the statement in Proposition 2.2(iii) to higher Picard number  we may need to use some 
more sophisticated number theory on algebraic varieties; for instance if $\rho =3$ we may use 
Siegel's Theorem concerning integral points on affine curves \cite{HS}.

Under the assumption that $X$ is general in moduli, in \it either \rm the case $c_2 (X) \cdot E <0$, \it or \rm the case $c_2 (X) \cdot E = 0$ and $E^3 >0$, we know 
all the possibilities 
 for the Type II or \Ttz exceptional surface $E'$ on the corresponding flopped model $X'$.  Moreover the numbers 
$n_d$ and the intersection multiplicities $E'\cdot \eta$ in the flopping formulae ($1'$) can take only  a small number of possible values.  By results of various  authors cited below, this 
enables us to classify locally analytically the possible flops which may be needed to pass from $E' \subset X'$ to $E\subset X$ and the way in  which flopping curves can 
intersect the surface corresponding to $E$.  We shall see that this in turn ensures that the extra codimension one singular 
locus introduced on $E$ via the flop has generically only finiitely many possibilities, and so for a generic hyperplane section $H$, the 
hyperplane section $C\subset H$ of $E$ has only finitely many possible types of singularity; hence for some fixed 
$\mu_0 >0$, the pair $(H, \mu C)$ is klt (Kawamata log terminal) for all (non-negative) real $\mu < \mu _0$.
It is this statement which is needed in Sections 6 and 7 for some boundary cases in the proof of Theorem 0.3,  when we wish to apply the 
Kawamata--Viehweg form of Kodaira Vanishing to a  real divisor on a generic 
hyperplane section of $X$, the fractional part of this divisor being supported on the given curve.  

In the following Lemma, we shall in fact be more explicit than is needed later.  The restriction to the numbers in 
the flopping formulae ($1'$) shows that the flops needed are locally just simple cases of the so-called Pagoda flops,
for which the geometry was explicity described by Reid in terms of blow-up and blow-downs, and 
which together with our knowledge of the possible  intersection numbers $E' \cdot \eta$ in the formulae ($1'$), gives 
the concrete fact that the codimension one singular locus consists generically of simple nodes or simple cusps.  Thus 
the above curves $C$ will only have at worst simple cusps or simple nodes as singularities, and in particular the reader will verify we may take  $\mu _0 =5/6$ in the above description.

\begin{lem}  We assume that the \CY threefold $X$is general in moduli.  For any very ample divisor class and 
divisor $H$ generic within that class, and with the rigid non-movable surface $E$ satisfying the numerical conditions {\bf either} $c_2 (X) \cdot E <0$, {\bf or}  $c_2 (X) \cdot E = 0$ and $E^3 >0$, the curve $C =E|_H$ on the smooth surface $H$ has at worst simple cusps or simple nodes as singularities. 
Thus in particular
the pair $(H, \mu C)$ is klt (Kawamata log terminal) for all (non-negative) real $\mu < {5\over 6}$. 
\begin{proof}
From our description above of the exceptional surface $E$ of a Type II or \Ttz contraction with $E^3 \ge 6$, we have at worst simple nodes on $C$  when $X=X'$, i.e. no flops are required.
Let us therefore consider first the case when a single flop is needed to go from $X'$ to $X$.
The crucial point is that  the formulae (${1'}$) and (3) then imply that $n_d =0$ for $d>2$.  If $n_2 \ne 0$, then 
$n_2 = 1$; $E'$ could not be a smooth scroll over $\P^1$ or a quadric cone, since then $c_2 (X)\cdot E \ge 0$ and $E^3 \le 0$.  Thus if $n_2 =1$, we have 
$E' = \P^2$ corresponds to a Type II contraction; if in addition $n_1 >0$ we would also deduce that $c_2(X)\cdot E \ge 0$ and $E^3 \le 0$, 
contrary to assumption.  Thus $n_1 =0$ and $n_2 =1$, and then as observed in Section 1 
 the flopping locus is  a single $(-1, -1)$-curve on $X$ (intersection multiplicity two with  $E'$), with the flop 
 being the Atiyah flop in this curve; it is then clear geometrically that 
 the generic hyperplane section of $E$ will have only simple nodes or  simple cusps (depending on whether the flopping curve intersects $E$ in two \it distinct \rm points or not).  We remark in passing that here and elsewhere in this proof, when the surface $E'$ is a (possibly singular) weak del Pezzo surface, i.e. with $-E'|{E'}$ nef, no component of the flopping curve can be contained in $E'$ (since by definition $E'$ is 
 strictly positive on the flopping locus).
  
If $n_2 =0$, then $n_1 \le 3$.  
If $n_1 =1$ then, again as observed in Section 1,  there is a unique flopping curve $l$ on $X'$, which is the $(-1, -1)$-curve in an Atiyah flop; it has intersection multiplicity at most two with $E'$, and in fact one unless $E' = \P^2$.  
When the flopping curve intersects $E' = \P^2$ twice, the  flop will as before introduce a  
codimension one singular locus on $E$, but the generic hyperplane section of $E$ will only have 
simple nodes or simple cusps.  When the flopping curve meets $E' =\P^2$ once, then $E$ is the blow-up of $E'$ in a point.
In the case when $(E')^3 =8$ and $l \not \subset E'$, both $E'$ and $E$ will be smooth in 
codimension one, as we are 
flopping in a single $(-1, -1)$-curve  meeting $E'$ transversely at a single smooth point.
In the case when the  surface $E'$ on $X'$ has $(E')^3 = 7$  (and therefore $c_2(X' )\cdot E' = -2$) and $l \not \subset E'$, it may be that we are in the case when $E'$ is non-normal, but we are then only allowed to flop in a single $(-1,-1)$-curve to obtain $X$, where the flopping curve has intersection multiplicity one with $E'$ 
and therefore does does not meet the double locus of $E'$ and only meets $E'$ in one point (transversely).  
Therefore the flop will introduce no further singularities to $E$.  From the description above of the possible non-normal surface $E$, the generic hyperplane section only has nodes as singularities.  We now  consider the cases when $l\subset E'$; from the remark at the end of the previous paragraph, here the only possibilities are \Ttz contractions, where $(E')^3 = 8$ or $7$.   Moreover the same remark ensures that when $(E')^3 =7$, 
we are not in the non-normal case, since if so we saw that $E'$ would be del Pezzo.   Therefore $E'$ is normal and we can use the classification of Lemma 3.1 
from \cite{WilGW} to deduce that $E'$ is either a Hirzebruch surface (when $(E')^3 =8$) or the blow-up of a Hirzebruch surface in one point (when $(E')^3 =7$).  The smooth surface contains $l \cong \P^1$ with 
$K_{E'} \cdot l  = E'\cdot l =1$, and hence we have $(l^2)|_{E'} =-3$; from this it follows that
  $E'$ is either the Hirzebruch surface $\F _3$ with $l$ the minimal section, 
or $\tilde \F_3$ (the blow up of $\F_3$ in a point not on its minimal section), corresponding 
respectively to $(E')^3 = 8$ or $7$, with $l$  the unique $(-3)$-curve on $E'$.  By Remark 5.13(a) in 
\cite{Reid}, the flop will introduce a curve of double points on $E$, 
which we claim consists generically of simple nodes.  For this last statement, 
observe that the flop is given by blowing up 
the $(-1, -1)$-curve $l$ to get an exceptional surface $F \cong \P^1 \times \P^1$, and that the intersection  $R$ of the proper transform of $E'$ with $F$ is a curve isomorphic to $\P^1$ which is a section with respect to the first ruling on $F$.  A straightforward calculation shows that it is a bisection with respect to the other ruling, and so under the flop, $E$ aquires a curve of double points which are generically simple nodes.  
   Thus in all cases 
where $n_1 \le 1$, the generic 
hyperplane section $H$ of $E$ will  have at worst either simple nodes or simple cusps as singularities.

We may assume therefore that $n_d =0$ for $d>1$ and 
$2 \le n_1 \le 3$.  If $n_1 =3$, then $E'$ must be $\P ^2$, but for $n_1 =2$ it can also be
 a smooth scroll over $ \P ^1$ or a quadric cone; 
 note that $E' \cdot \eta =1$ in the flopping formula ($1'$).  
Let $Z$ denote the exceptional locus of the Type I contraction on $X'$ corresponding to the flop. When $E' = \P^2$, no component of $Z$ lies in $E'$.  
Assuming $Z$ to be connected,
we saw in Section 1 that the image $\bar U$ of an appropriate neighbourhood $U$ of $Z$ is then a compound Du Val singularity.

We assume first that $Z$ is just a single $\P^1$.  One possibility is  that 
the generic  section $Y_0$ through the singularity  $\bar U$
has just an $A_1$-singularity with resolution $X'_0$, and that the map $f : \Delta \to \hbox{Def}\; X'_0$ has ramification index $m$.  
Here we are in the situation of (1.4) from \cite{Reid}; we note that  $U \to \bar U$ is a simultaneous resolution of the family $\bar U$, and we know that 
${\rm Def} (X_0 ') \to {\rm Def} (Y_0 )$ is a ramified double cover, corresponding to the Weyl group $C_2$ of the singularity $A_1$.  Moreover, the map $f : \Delta 
\to {\rm Def} (X'_0 )$ may be written in the form $f(t) = t ^m$, the map $g : \Delta 
\to {\rm Def} (Y_0 )$ in the form $g(t) = t^{2m}$ for some local analytic coordinate $t$ on $\Delta$, 
and $\bar U$ may be written locally analytically as the hypersurface singularity $x^2 + y^2 + z^2 + t^{2m} =0$ (cf. Corollary 1.16 of \cite{Reid}).  Deforming $f$ so that it is tranverse to the origin, we see that $m =n_1$.  In the language of Section 5 from \cite{Reid}, this says that $Z \cong \P^1$ 
is a $(-2)$-curve which has  \it width \rm $n_1$ (see Definition 5.3), and the analytic neighbourhood of $Z\subset X'$ is determined by this data (Remark 5.13(b)).  
Moreover the corresponding flop (a \it Pagoda \rm flop) is explicitly described there as the composite of $n_1 -1$ blow-ups in 
$(-2,0)$-curves, an Atiyah flop in a $(-1, -1)$-curve, followed by $n_1 -1$ blow-downs to $(-2,0)$-curves (see Pagoda 5.8 in \cite{Reid}).  Since $n_1 =2$ or $3$, the flop is therefore analytically one of the
two simplest flops in a $(-2,0)$-curve; if the width were one we would have an Atiyah flop in a $(-1,-1)$-curve.  In particular, since $E'$ has intersection multiplicity one with
the reduced curve $Z\cong \P^1$, we know that \it  either \rm
$Z$ intersects $E'$ transversely in one (smooth) point \it or \rm  $Z\subset E'$.  In the former case, 
from the given geometric description of the flop from \cite{Reid}, we see by an argument analogous to that for an Atiyah flop in a $(-1,-1)$-curve, that 
the only new singularity introduced on $E$ via the flop is an $A_1$ 
or $A_2$ rational double point (according to whether $n_1 =2$ or $3$).
 In particular $E$ has no codimension one singular locus.  In the  case when $Z\subset E'$, a previous argument ensures that
 we have  a \Ttz contraction, for which $(E')^3 =8$ and hence $n_1 =2$ 
 (so that analytically the flop is the simplest flop in a $(-2,0)$-curve).  In particular the exceptional locus is a Hirzebruch surface; since $K_{E'}\cdot Z = E'\cdot l =1$, we deduce 
$E' \cong \F_3$ with $Z^2 =-3$ on $E'$ (compare with the previous argument with $n_1 =1$).  From  Remark 5.13(a) of \cite{Reid}, the flop will introduce a curve of double points on $E$; more precisely, when 
we blow up $Z$, we obtain an exceptional surface $F \cong \F_2$, whose minimal section $M$ is a $(-1,-1)$-curve 
on the blown-up threefold; a straightforward check verifies that 
the intersection of the proper transform of $E'$ with $F$ is a section (isomorphic to $\P^1$) of $F$ which has intersection multiplicity one with $M$.  When we make the Atiyah flop in $M$, the proper transform of $E'$ 
therefore now contains the minimal section of the resulting $\F_2$ as well as a second (distinct) section arising from the first blow-up of $Z$.  Thus contracting the $\F_2$ to complete the factorisation of the Pagoda flop, shows that the double locus on $E$ consists generically of  simple nodes.

When the exceptional locus of the Type I contraction is a $\P^1$, a second possibility might be $Y_0$ 
having a more complicated
singularity (and so $ X'_0$ is a partial resolution), for which however there is a finite classification \cite{KatzM, KawHyp}.  If 
$Y_0$ is not an $A_1$-singularity, it will be one of $D_4 , E_6. E_7 , E_8$, where in the first three cases the node 
of the Dynkin diagram corresponding to the non-contracted 
curve is uniquely specified, and for $E_8$ there are two possibilities for the node \cite{KatzM}.  
Using the description from Section 2 of \cite{KatzM} of the root systems and the action of $W_0$ 
on the roots, we see that the number of irreducible components $r$ 
of the discriminant locus $D$ on $\hbox{\rm Def}\, X'_0$ has $r>3$ in all these other cases, and hence in particular $\delta >3$, and so they do not occur 
in our case (where $n_1 \le 3$).

Continuing the case of a single flop with $n_d =0$ for $d>1$ and 
$2 \le n_1 \le 3$, we might have the flopping locus being disconnected.  Here
we just locally argue separately on each connected component.  Each component of the flopping locus will 
correspond of one or two $(-1, -1)$-curves when we deform the complex structure locally in a neighbourhood of that component, from which we deduce 
 that the components of the exceptional locus are isomorphic to $\P^1$
(compare next paragraph, which shows there would otherwise 
be at least three $(-1, -1)$-curves when we deform the complex structure locally),
and that the local flops must then  either be Atiyah flops  or Pagoda flops as described above with $m=2$; from Reid's geometric 
description of these flops we again conclude that $E$ has at worst a curve of singularities which are generically nodes.

The only remaining case with a single flop is when 
$Z$ is connected but has two or more irreducible components.  In this case 
we can calculate that the discriminant divisor $D$ on $\hbox{\rm Def}\, X'_0$ 
has at least four components, with the only exception being the case of a full resolution of an $A_2$-singularity, and hence $Z = C_1 \cup C_2$, 
when $D$ has three components.  Since $n_1 \le 3$ and hence $r\le 3$, we must be in this latter case and our assumptions 
then imply that  $n_1 = 3$ and $f$ is unramified and generic.  We note also $(E')^3 =9$ and so $E' \cong \P^2$.  Moreover a small deformation of the complex structure on $U$ yields three $(-1,-1)$-curves $\Gamma _i$ ($i=1,2,3$), where the numerical assumptions of the Lemma imply that $E' \cdot \Gamma _i =1$ for all $i$.  We claim that this is impossible and so our assumptions imply that if $Z$ is connected, then it must be irreducible.

The above claim follows from the calculations in Section 4 of \cite{KatzM}, and for ease of explanation we adopt notation and conventions which are close to those in that paper.  
Here ${\rm Def} (X_0)$ is the subspace of $\C ^3$ defined by the linear form 
$t_1 + t_2 + t_3 =0$, over which we have a 
family of surfaces $ \calY$ with equation $$-xy + (z +t_1)(z+t_2)(z+t_3),$$ where the simultaneous resolution $\calZ \to {\rm Def} (X_0)$ may be constructed by taking 
the closure of the 
graph of a certain morphism $\calY \to \P^1 \times \P^1$, given explicitly in \cite{KatzM}.  Moreover it is shown there  that the three lines $t_i =t_j$ ($i\ne j$) of the discriminant locus 
correspond to the loci where $C_1$, $C_2$ and $C_1 \cup C_2$ deform to curves isomorphic to $\P^1$.  Therefore when we deform $f : \Delta \to {\rm Def} (X_0)$, and hence 
the complex structure on $U$, we obtain curves $\Gamma _i$, where we may assume that as the deformation parameter tends to zero, we have $\Gamma _1$ specializes to $C_1$, 
$\Gamma _2$ specializes to $C_2$ and $\Gamma _3$ specializes to $C_1 \cup C_2$; in particular it follows that $E' \cdot \Gamma _3 >1$, contrary to assumption.

There may of course be more than one flop needed to get from $X'$ to $X$, but those cases will be simpler than when there is only one flop.  Again we see that $E'$ is either $\P^2$, a smooth scroll over $\P^1$, or a quadric cone.
If we have a sequence of three flops, then $E'\cong \P^2$,  
and the formulae ($1'$) imply that each will be an Atiyah flop in some $(-1, -1)$-curve with intersection multiplicity one with the relevant transform of $E'$; no singularities will be introduced on $E$.  If we have a sequence of two flops, then one of them might have 
$n_1 =2$ and the other $n_1 =1$, or both might have $n_1 =1$.  From the Reid description of the Pagoda flop with $n_1 =2$, $E$ will only have a curve of singularities if $E'$ is the ruled surface $\F_3$ or $\F_2$, and the Reid flopping curve is a $(-3)$-curve on the corresponding surface $\F_3$ or $\tilde \F_3$ defined before.  In this case also the generic singularity on the singular locus of $E$ is a simple node.

 In conclusion,  for any very ample divisor class and 
divisor $H$ generic within that class, 
the curve $C =E|_H$
on the smooth surface $H$ has only simple nodes or simple cusps as singularities,  
although the number of such singularities depends on the divisor class of $H$.  
Thus in particular 
the pair $(H, \mu C)$ is klt  for all $\mu < {5\over 6}$.  
 \end{proof}  \end{lem}
 
  We prove now a further inequality on rigid non-movable surfaces, which will be needed when dealing with a special case in Theorem 6.2.
 
 \begin{lem} Let $E$ be a rigid non-movable surface on a \CY threefold $X$. then 
$$ E^3 + {1\over 2}(c_2(X)\cdot E)^3 \ge - 99.$$
 \begin{proof}
 We claim first that if $E_1$ is a rigid non-movable surface on a \CY threefold $X_1$, and $X_2$ is obtained from $X_1$ by an $E_1$-directed flop, then 
 $$ E_1^3 + {1\over 2}(c_2(X_1)\cdot E_1)^3 \ge E_2^3 +  {1\over 2}(c_2(X_2)\cdot E_2)^3.$$
 Note from ($1'$) that $E_1^3 = E_2^3 - (E_2 \cdot \eta)^3 \sum_{d>0}n_d d^3$ and 
 $c_2(X_1)\cdot E_1 = c_2 (X_2)\cdot E_2 + 2(E_2 \cdot \eta) \sum_{d>0}n_d d$, where $E_2\cdot \eta >0$.  If we set $N = 
 (E_2\cdot \eta) \sum_{d>0}n_d d$, then 
 $$ E_1 ^3 \ge E_2 ^3 -N^3, \quad c_2(X_1)\cdot E_1 = c_2(X_2)\cdot E_2 + 2N.$$
 Therefore 
 $$  (c_2 (X_1)\cdot E_1 )^3 =  (c_2 (X_2)\cdot E_2 )^3 +6N (c_2 (X_2)\cdot E_2 )^2 + 12 N^2 
 (c_2 (X_2)\cdot E_2 ) + 8N^3.$$ Hence 
 $$ E_1^3 + {1\over 2}(c_2(X_1)\cdot E_1)^3 \ge E_2^3 + {1\over 2} (c_2(X_2)\cdot E_2)^3 + NG,$$
 where $$G = 3N^2 + 6 N (c_2(X_2)\cdot E_2)  + 3 (c_2(X_2)\cdot E_2)^2 = 3(N + c_2(X_2)\cdot E_2)^2 \ge 0.$$
 
 Given now $E$ as in the statement of the Lemma, we can make a sequence of directed flops to achieve a 
 surface $E'$ on a birational model $X'$,  where $E'$ is the exceptional surface for a Type II or \Ttz contraction.
 From equations (2) and (3), we check easily that the minimum value of $(E')^3 + {1\over 2}(c_2(X')\cdot E')^3$ 
 occurs when $c_2(X')\cdot E' = -6$, and so $(E')^3 + {1\over 2}(c_2(X')\cdot E')^3 \ge 9-108 =-99$; the lemma now follows by repeated use of the previous claim.
\end{proof}  \end{lem}
 
A final question concerns how many rigid non-movable surfaces on $X$ there can be.  An easy argument shows that for any given $\rho +1$ irreducible surfaces on $X$, 
some integral combination of at most $[(\rho +1)/2]$ of them must be mobile, where $[\, .\, ]$ denotes the integral part of a number  (there will be an integral dependence between the corresponding classes, and then just take terms 
with negative coefficients to the other side).  Remark 1.2 implies that  
$c_2 (X)$ is non-negative on any mobile class on $X$; this implies that we cannot have more than $\rho$ rigid non-movable surfaces $E$ with 
$c_2 \cdot E <0$.  
In fact we cannot have more than $\rho$ such surfaces $E$ with $c_2 \cdot E \le 0$, unless there is a base point free linear system $|D|$ on $X$ 
with $c_2 (X)\cdot D = 0$; $X$ then has a very special 
structure, as studied in \cite{Og}.  
In general there may exist infinitely many rigid non-movable surfaces, as for instance would be the case with an elliptic \CY threefold over a rational surface $S$ for 
which the generic fibre over $\C (S)$ had infinite Mordell--Weil group.  A natural question to ask if whether any \CY threefold $X$ containing  infinitely many rigid non-movable surfaces 
has to be of fibre type.  In the case we study in detail later in this paper, namely $\rho =2$, it will  however be clear that there are at most two rigid non-movable surfaces.

\section{Components of the positive index cone}

In this section we assume that $X$ is a smooth complex projective threefold and that  the cubic and linear forms on $H^2 (X, \Z)$ are given.  When $X$ is a \CY threefold, 
 to obtain boundedness, we are trying to find 
 a specific set of rational classes in $H^2 (X, \Q)$, at least one of which will be movable and big on 
 $X$.  We shall furthermore ignore the easy case when $\rho (X)=1$.  
We then define two more open cones in $H^2 (X, \R)$.  For $D\in H^2 (X, \R)$, we denote by $q_D$ the 
associated quadratic form on $H^2 (X, \R)$ given by $q_D (L) = D\cdot L^2$.  Saying that the quadratic form $q_D$ is singular is just saying that the Hessian 
determinant of the cubic form vanishes at $D$.

We define the \it index cone \rm to consist of elements $D \in H^2 (X, \R)$ for which the quadratic form $q_D$ has index $(1, \rho -1 )$.  The \it positive index cone \rm is defined to be the subcone of the index cone on which the cubic form is positive.

\begin{lem}
For a homogeneous form $G$ of degree $d$ on $H^2 (X, \R) = \R ^\rho$, the complement of the hypersurface defined by $G$ consists of a finite number $n$ of connected components, where there an explicit bound on $n$ in terms of $d$.

\begin{proof} 
This is a direct consequence from Theorem 2 of \cite{Warren}.
\end{proof}\end{lem}

 We can apply Lemma 3.1 in the case when $G$ is the Hessian form for the cubic.  The index cone is then a finite (disjoint) union of open cones corresponding to a certain subset of the connected components of the complement of the Hessian hypersurface.  We can also apply Lemma 3.1 in the 
 case when $G$ is the product of the cubic form and its Hessian form.  
 The positive index cone will then 
 be a finite (disjoint) union of open cones corresponding to a certain subset of the connected components of the complement of the hypersurface $G$.
 
 This cone will be important in the second part of this paper dealing with the \CY case with $\rho =2$, but more indicative of the general \CY case is $\rho =3$, which is studied in detail in \cite{WilCY3}.  If for instance the 
 cubic is smooth and the corresponding real elliptic curve in $\P ^2 (\R )$ has two connected components, 
 there are four components of the positive index cone, one corresponding to the positive cone on the bounded 
 component of the real elliptic curve and three \it hybrid \rm components, part of whose boundary is given by 
 the cubic form and part of whose boundary is given by the Hessian form.  
 
 For $X$ a smooth complex projective threefold with a given cubic form, there is a relationship between the cones just defined and cones introduced in Section 1.
 
 \begin{lem}
 If $X$ is a smooth complex projective  threefold with a  given cubic form on $H^2 (X, \R)$, then any class 
 $D \in \Mov (X)$ not lying on the Hessian hypersurface is contained in the index cone.  Moreover, 
  the \Kahler cone $\calK$ is contained in a connected component of the positive index cone.
 \begin{proof}
 We know that any big mobile divisor $D$ may be taken to be irreducible by Bertini's theorem  and then it may be easily checked (by pulling back divisors to a desingularization and using the Hodge index theorem) that the associated quadratic form $q_D$ has index $(1,t)$ for some $t \le \rho -1$.  Thus 
 unless the class $D$ lies on the Hessian hypersurface to the cubic form, the index of $q_D$ will be $(1, \rho -1 )$.  Therefore a big movable rational class $D$ will have associated index $(1, \rho -1 )$ unless it lies on the Hessian hypersurface, since some integral multiple will be mobile and big.  Since however the index 
cone is open, the statement remains true for all real classes $D \in \Mov (X)$.
 
 By the Hard Lefshetz theorem, the Hessian doesn't vanish at any class in $\calK$, and thus $\calK$ is contained in a connected component of the index cone.  Since the cubic form is also positive on the convex 
 cone $\calK$, it follows that $\calK$ is contained in a component of the positive index cone.
  \end{proof}
 \end{lem}
 
 In the case when the cubic comes from the cup-product on $H^2 (X, \R) = \R ^3$ for $X$ a \CY threefold with $\rho =3$, the manuscript \cite{WilCY3} studies which of these components of the positive index cone might be 
 expected to contain the \Kahler cone.
 
   Let 
us now fix a connected component $P^\circ$ of the positive index cone ---  there are 
 only finitely many possibilities for $P^\circ$ (one of which will contain the \Kahler cone).
We let $P$ denote the closure of $P^\circ$.

\begin{lem} If $D_1 , D_2 , D_3 \in P^\circ$, then $D_1\cdot D_2\cdot D_3 > 0$.  In particular, a convex combination of two classes in $P^\circ$ has strictly positive cube.

\begin{proof} We show first that $D_1 ^2 \cdot D_2 >0$ for all $D_1 , D_2  \in P^\circ$; suppose not, so that  there are (linearly independent) $D_1, D_2 \in P^\circ$ with $D_1  ^2 \cdot D_2 \le 0$.  We can join $D_2$ to $D_1$ by a real curve (say piecewise linear) in $P^\circ$, on which there are only finitely many points $D$ with $D_1 ^2 \cdot D = 0$.  Suppose that $D$ is the \it last \rm point (going from $D_2$ to $D_1$) for which this occurs
(possibly including the case $D=D_2$).  Since $D_1^3 >0$, the index condition at $D_1$ 
implies that $ D^2\cdot D_1 < 0$.  Furthermore we deduce then 
that there exists $D'$ on the curve strictly between $D$ and $D_1$  for which $D'^2 \cdot D_1 =0$. Since $(D')^3 >0$, the index condition at $D'$  implies that $D' \cdot D_1 ^2 <0$.  Thus for some $D''$ strictly between $D'$ and $D_1$ on the curve, we have $D'' \cdot D_1 ^2 =0$, contradicting 
our initial assumption on the choice of $D$.

Thus $D_1^2 \cdot D_2 >0$ for all $D_1 , D_2 \in P^\circ$; 
now suppose $D_1 \cdot D_2 \cdot D_3 \le 0$ 
and join $D_3$  to $D_1$ by a curve, and let $D$ denote a point 
on the curve  for which $D_1 \cdot D_2  \cdot D = 0$.  Since $D_1 ^2 \cdot D >0$, the index 
condition at $D$ implies that $D_2^2 \cdot D <0$, contradicting what we proved in the earlier part.
  \end{proof}\end{lem}

The \Kahler cone $\calK$ of $X$ will be contained in one of these connected cones $P^\circ$.  However for any surface $E$ on $X$, we also have $D^2 \cdot E >0$ for all $D \in \calK$.  In what 
follows below and in subsequent sections, we shall take $E_1 , \ldots E_r$ to be the classes of a rigid non-movable surfaces; here we prove a more general result.

\begin{prop} Suppose $E_1, \ldots , E_r$ are classes such that the 
associated quadratic form to each $E_i$ has index $(1, s_i)$ 
for some $s_i \le \rho -1$ 
(for instance the classes of irreducible surfaces).  For each of the connected components $P^\circ$ 
of the positive index cone, 
 the extra conditions that $D^2 \cdot E_i  >0$ for $i = 1, \ldots , r$
yield only a finite number of connected open subcones of $P^\circ$;
  for $D_1 , D_2$ in such a subcone, we have $D_1 \cdot D_2 \cdot E_i  >0$ for all 
   $i = 1, \ldots , r$.

\begin{proof}
We follow the same ideas as in the proof of Lemma 3.3.  The condition $E_i \cdot D^2 =0$ 
defines a further hypersurface, defined by some quadratic form.  We consider the 
complement of the 
 hypersurface in $\R^\rho$ defined by the product 
of these $r$ quadratic forms, the
cubic form and its Hessian.  By Lemma 3.1 this is an open cone with finitely many connected components, 
where again there is an explicit upper bound on the number of such connected components.    This then subdivides  $P^\circ$ into a finite number of open connected cones, on each of which $D^2 \cdot E_i \ne 0$ for all $i$.  
In particular this yields a finite number of connected open subcones on which $D^2 \cdot E_i >0$
for all $i$.

Suppose $D_1, D_2$ are in such a subcone $Q$ but $D_1 \cdot D_2 \cdot E_i \le 0$
for some $i$.  
Joining $D_1$ to $D_2$ by a real curve in $Q$, we may find $D$ on this curve with 
$D_1 \cdot D \cdot E_i = 0$.  Since by construction $D_1 ^2 \cdot E_i >0$, 
the index assumption on the associated quadratic form to $E_i$ implies that 
$D ^2 \cdot E_i \le 0$; this is the required  contradiction.
\end{proof} \end{prop}

Assuming that $P^\circ $ is the connected component containing the \Kahler cone  $\calK$, 
and $E_1 , \ldots , E_r$  are classes of rigid non-movable surfaces on $X$, there are only finitely many connected open subcones $Q$ of $P^\circ$ defined by the extra conditions 
$D^2 \cdot E_i >0$ for $i = 1, \ldots r$.  The \Kahler cone will be contained in one of these subcones.  In particular, if $\calK$ is contained in such a subcone $Q$, and we take 
$H$ to be a very ample divisor class, then $D\cdot H \cdot E_i >0$ for all $D \in Q$ and 
$i = 1, \ldots , r$.

In the case when $X$ is  a \CY threefold which is general in moduli and $\{ E_1 , \ldots , E_r \}$ is a complete set of rigid non-movable classes, we observe in the next section  for $H$ a generic very ample divisor that $D|_H$ is nef for all $D\in Q$.  If $D$ is any rational class in $Q$, then we shall apply cohomological methods below to obtain effectivity of $mD$ (and indeed 
a non-trivial mobile part) for some integer $m$ depending only 
on $D^3$ and $c_2(X)\cdot D$.

\section{Finding movable classes}

With the notation as in Section 3, we shall assume 
in this section that $X$ is a \CY threefold with $\rho >1$ and that the component 
$P^\circ$ of the positive index cone 
contains the \Kahler cone; there are only finitely many possible  components $P^\circ$ of the positive index cone, so 
for the purposes of proving boundedness this 
assumption may be made without loss of generality.

\begin{prop} Suppose $H$ is a very ample smooth divisor on 
a \CY threefold $X$ and we are given an integral class $D \in P^\circ$, the component of the positive index cone containing the \Kahler cone.  If $D|_H$ is nef, then $h^0(X, \O_X (mD))  >1$ for some integer $m>0$, depending only on 
$D^3$ and $c_2(X)\cdot D$,
and moreover $D$ is big.

\begin{proof} This is an argument using the Riemann--Roch  and vanishing theorems.  We note 
from Lemma 3.3 that $D^2\cdot H >0$, $D\cdot H^2 >0$, and so the Riemann--Roch theorem on $H$ implies that 
$$h^0 (H, \O_H (mD))  \ge \frac{1}{2} mD\cdot(mD -H)\cdot H + \chi (\O_H)  \sim \frac{1}{2}m^2 D^2\cdot H$$ 
for $m\gg 0$.  Thus $D|_H$ is big. 
Using  also the nefness 
of $D|_H$ however, we deduce that $h^1 (H, \O _H (rH + D)) = 0$ for all $r\ge 1$.  Now take the long exact sequence of cohomology associated (for $r\ge 1$) to the short exact sequence of sheaves
$$ 0\to \O _X (D + (r-1)H) \to \O _X (D + rH) \to \O _H (D + rH) \to 0$$
to deduce that $h^2 (X, \O _X (D + (r-1)H)) \le h^2 (X, \O _X (D + rH))$ for all $r\ge 1$.  Serre vanishing implies that $h^2 (X, \O _X (D + rH)) =0$ for $r \gg 0$, and so $h^2 (X, \O _X (D)) =0$;  similarly $h^2 (X, \O _X (nD)) = 0$ for all $n \ge 1$.

Now use the Riemann--Roch theorem on $X$ to deduce that $$
h^0(X, \O _X (mD)) \ge \chi (\O _X(mD)) = {1\over 6} m^3 D^3 + {1\over 12} m D \cdot c_2 (X),$$ 
which is $>1$ for some positive integer $m$ depending only on $D^3$ and $D\cdot c_2$;
 moreover $h^0(X, \O _X (nD)) \sim n^3$ for large $n$, i.e. $D$ is big.\end{proof}\end{prop}

\begin{rem}
Given an integral divisor $D$ on $X$ and a positive integer $m$ such that 
$h^0(X, \O_X (mD))  >1$, we can write 
write $| mD | = | \Delta | + \calE$, where  $\Delta$  is a non-trivial effective divisor 
  not having any rigid non-movable surface as an irreducible  
component, whilst $\calE$ is supported on the set of rigid non-movable surfaces.  
Assuming  $X$ to be general in moduli, any irreducible component of 
$\Delta $ will correspond on some 
minimal model to a nef divisor (and hence semi-ample by \cite{Og}),  and thus some multiple of it is a mobile class on $X$.  Hence some multiple of $\Delta$ is also mobile and 
in particular $\Delta$ is movable.
The idea now is to work with this movable  divisor $\Delta$, despite the fact that we don't know that $\Delta$ itself is necessarily big in general.  The proof of Theorem 4.5 below illustrates the nature of the argument that will be employed in subsequent sections.
\end{rem}

 For $H$ generic in its linear system, and hence smooth by Bertini's theorem, we  relate the condition that $D|_H$ is nef to the rigid 
non-movable surfaces.

\begin{lem}  Suppose that a \CY threefold $X$ is general in moduli and there are only finitely many rigid 
non-movable surfaces $E$ on $X$.  Moreover let $Q$ denote the 
particular connected component of the open subcone of $P^\circ$ consisting of classes $L$ with $L^2 \cdot E >0$ for all the rigid non-movable surfaces $E$, with $Q$ also containing the \Kahler cone.  Let $|H|$ denote a very ample linear system on $X$.  If $D \in \bar Q \cap P^\circ$, then 
 $D|_H$ is nef for a generic element $H$ of the linear system.

\begin{proof}  It follows from the proof of Proposition 4.1 that when $H$ is smooth, 
the restriction $D'|_H$ is big 
for all rational classes $D'$ sufficiently close to $D$ (so that $D' \in P^\circ$); thus by definition we have 
$D|_H$ is big for $H$ generic in its linear system. 
 In particular it may be written as a finite sum of curves with 
positive real coefficients, and it is only on these curves that $D|_H$ could be negative.  
If $D|_H$ were negative on some curve on $H$, then our assumption on $H$ being generic in its linear system ensures that by varying $H$ in the linear system, a surface $E$ is swept out by such a curve.  Given that $X$ is assumed to be general in moduli, it follows that $E$ must be 
a rigid non-movable surface.  
By Proposition 3.4, we know however that $D|_H \cdot E|_H = D\cdot H \cdot E \ge 0$ 
for all the rigid non-movable surfaces $E$, and  hence 
 $D|_H$ is nef.
\end{proof}\end{lem}

A similar proof to Proposition 4.1 yields a fact about rigid non-movable surfaces.

\begin{prop}  Suppose $P^\circ$ is the component of the positive index cone which contains the \Kahler cone (with corresponding closed cone $P$) and $E$ is the class of a 
rigid non-movable surface on $X$; then $E \not \in P$ 
unless both $E^3 =0$ and $c_2 \cdot E =0$.

\begin{proof}  Suppose $E \in P$.  From Lemma 3.3 it follows that $E^2 \cdot H \ge 0$ for any very ample class $H$.  Thus 
$E|_H$ is nef on any very ample smooth divisor $H$, generic in its linear system. 

If the linear form $E^2$ is trivial, i.e. 
$E^2 \equiv 0$, then
$E$ will be nef, since the prime divisor $E$ is non-negative on any curve not in $E$, and it is zero on any curve in $E$.    We can then deduce from \cite{Og} that 
some positive multiple of the effective class $E$ is mobile (in fact free), a contradiction.

We can therefore choose the very ample class $[H]$ so that $E^2 \cdot H > 0$.
With $H$ then chosen generic in its linear system, we have $H$  smooth and irreducible by Bertini's theorem; 
moreover $E|_H$ is nef and big.
For $n>0$, we use the short exact sequence of sheaves 
$$ 0 \to \O _X (nE) \to \O _X (nE +H) \to \O_H (H + nE) \to 0.$$
Therefore $h^2 (H, \O_H (H + nE)) = h^0 (H, \O_H (- nE))  =0$ and the Kodaira Vanishing theorem implies that $h^1 (H, \O_H (H + nE)) =0$; hence $h^2(X, \O_X (nE)) = h^2 (X, \O_X(nE +H))$; by taking a large enough multiple of the class $[H]$, we may without loss of generality assume that 
 this latter term is zero.  Thus the Riemann--Roch theorem yields  that $$h^0(X, \O_X (nE)) \ge \chi (X, \O_X (nE)) = {1\over 6} n^3 E^3 + {1\over {12}}n c_2 \cdot E.$$

We now have two cases.  If $E^3 >0$, we can choose  $n>0$ such that $2n^3 E^3 + nE\cdot c_2 \ge 24$, from which it follows that $h^0(X, \O_X (nE)) >1$, contradicting 
non-movability.  We must therefore have $E^3 =0$, and then from Proposition 2.2 (i)  that $c_2 \cdot E \ge 0$.  If however $E^3 =0$ and 
$c_2 \cdot E >0$, then we can still choose $n$ so that 
$2n^3 E^3 + nE\cdot c_2 \ge 24$, from which it follows that $h^0(X, \O_X (nE)) >1$, contradicting non-movability.  Thus we are left only with the case $E^3 =0$ and $c_2 \cdot E =0$.
\end{proof}\end{prop}

\begin{thm} Suppose that in addition to the cubic and linear forms on $H^2 (X, \Z)$ we are given $\rho +1$ rigid non-movable surface classes $E$; 
then {\bf either} some birationally equivalent minimal model 
$X'$ is of fibre type  {\bf or }$X$ lies in a bounded family.  If there are at most $\rho$ rigid non-movable surfaces on $X$,  all of whose classes are specified and have the property 
that each class $E$ 
satisfies $c_2 \cdot E >0$,   then (assuming $X$ to be general in moduli) the same conclusion holds.

\begin{proof} 
When there are  $\rho +1$ classes $E$, they are linearly dependent, which then shows that some positive convex combination 
$\Delta$ of these divisors moves (cf. argument from the end of Section 2).  
Provided we have started from  a minimal linear dependence relation, the linear system $|\Delta | $ does not have any fixed component, and so $\Delta$ is mobile.  There are only finitely many such classes $\Delta$ arising in this way, and we just choose one of them.

  If our chosen $\Delta$ is also big, then we have boundedness by Proposition 1.1.  If 
$\Delta$ is not big, then we are in the fibre-type case (a component of a generic element of $| \Delta |$ gives rise to a nef, but not big, irreducible divisor 
on some smooth minimal model $X'$ birationally equivalent to $X$, and hence by \cite{Og} to a fibre space structure on $X'$).

We assume therefore that there are at most $\rho$ classes $E$ corresponding to rigid non-movable surfaces and that these classes are also specified, 
and each such class $E$  
satisfies $c_2 \cdot E >0$.  From Lemma 3.1 there are only finitely many possibilities for the cone $P^\circ$ containing the K\"ahler cone, 
 and by Proposition 3.4 such a cone $P^\circ$ gives rise to finitely many open subcones $Q$ 
satisfying the conditions that $D^2 \cdot E > 0$ for all the (given) rigid non-movable surface classes $E$ on $X$.  

For any such cone $Q$, we specify a 
fixed choice of an integral divisor $D$ in $Q$.  As there are only finitely many cones $Q$, we have 
a finite set of (known) integral divisor classes.  We may assume without loss of generality that a given cone $Q$ contains the K\"ahler cone, and 
that $D$ is our chosen integral  class in $Q$.  Moreover we are assuming that $X$ is general in moduli.

Let $H$ be a generic very ample divisor on $X$; we know from Lemma 4.3 that  $D|_H$ is nef and big.  We now apply Proposition 4.1 to produce an $m>0$, depending only on 
$D^3$ and $c_2(X)\cdot D$, 
 with $h^0 (X, \O_X (mD)) >1$.  Under the assumption that $X$ is general in moduli, by 
 Remark 4.2 we can write $| mD | = | \Delta | + \calE$, with $\Delta$ 
 non-trivial movable with some multiple mobile 
 and $\calE$ is  supported on the (finite) set of rigid non-movable surfaces.  

We observe that $c_2 \cdot \Delta \ge 0$ by Remark 1.2,  
and that the assumptions of the theorem include that $c_2 \cdot E 
> 0$ for each of the (known) rigid non-movable surfaces $E$.  Moreover both $D$ and $m$ are 
known, and hence we also know 
$mD\cdot c_2$; this ensures that there are only finitely many possibilities for the coefficients of $\calE$, and hence for the fixed divisor $\calE$.  

Thus there are also only finitely many possibilities for the class of the movable part $\Delta$.  If now $\Delta$ is also big, the fact that we know the values $\Delta ^3$ and $c_2 \cdot\Delta$ up to 
finitely many possibilities shows via 
Proposition 1.1 that $X$ lies in a bounded family.  (Observe that in the case when there are no rigid non-movable surfaces on $X$, we are trivially in this case since $\Delta = mD$ is
movable and  big.)

In the exceptional case where the movable divisor $\Delta$ is not big, this will also be true for any of its components; such a component will correspond to a nef but not big divisor on some birationally equivalent minimal model $X'$, where  $X'$ will be of fibre type by \cite{Og}.\end{proof}\end{thm}

\begin{rem}
If  $X'$ is of fibre 
type, then it is either elliptic, a K3 fibre space  over ${\bf P}^1$ or an abelian fibre space  over ${\bf P}^1$.  For the first of these fibre types, we have at least birational boundedness from 
\cite{GrEll}; for the fibre spaces over ${\bf P}^1$ we in general do not even know this.  Given the extra information in the hypotheses of the theorem, we would however still hope to prove boundedness (as we do in the case $\rho =2$ later), 
but we do not pursue the question at this stage.
We remark in respect of the second case in the theorem that rigid non-movable surfaces $E$ with $c_2 \cdot E \le 0$ will be very special, 
as can be seen from the arguments of Section 2 (cf. Proposition 2.2).  \end{rem}
\
\noindent \it Proof of Theorem 0.1. \rm
When there are no rigid non-movable surfaces on $X$, we saw in the proof of  
the second part of Theorem 4.5 that we can choose a finite set of integral divisor classes $D$ (one for 
each of the connected components $P^\circ$ of the positive index cone), 
one of which will be movable and big (see remark in brackets at the end of the penultimate paragraph of the proof of Theorem 4.5) for $X$ general in moduli.  Applying Proposition 1.1,  we have 
therefore now proved 
Theorem 0.1  as stated in the Introduction. \hfill $\Box$
\smallskip

There is an  alternative proof of Theorem 0.1, which we shall now also give since the technique employed will be useful later, and it will moreover yield a stronger version of Corollary 0.2.  
We use the concept of the \it volume \rm $\vol (D)$ of 
a real divisor $D$, as explained in \cite{Nak} or Section 2.2 of \cite{Laz}, Vol 1.  We recall that the volume only depends on the numerical class of a divisor, and that a divisor being big is equivalent to its having strictly positive volume.  Moreover, as a function on the space of real numerical divisor classes, the volume is a continuous function (\cite{Laz}, Vol 1, Theorem 2.2.44).
\par
\

\noindent \it 2nd proof of Theorem 0.1. \rm We show that given the connected component $P^\circ$ of the positive index cone containing the K\"ahler cone, there can exist 
only finitely many families.  To see this, we first show that under the assumptions of Theorem 0.1, we have $\Movb (X) \supset P$.  For suppose not, then some part of the boundary of $\Movb (X)$ is in $P ^\circ$.  The formula (1) for how $L^3$ changes under flops implies that for any 
strictly movable class $L$, we have the inequality $\vol (L) \ge L^3$ --- since $\vol$ is continuous, this in fact is true also for $L\in \Movb (X)$.  Moreover, the continuity of $\vol$ implies 
that there exists a class $M \in P^\circ \setminus \Movb (X)$ such that $\vol (M) >0$.  Using the $\sigma$-decomposition of the big divisor $M$
(see \cite{Nak}, Chapter 3) in terms of a sum of positive and negative parts, the negative part is effective and non-zero (if zero then $M$ would be movable).   Assuming that $X$ is 
general in moduli, the components of the negative part of $M$ must be rigid non-movable surfaces, contradicting the assumption in the theorem.

Thus we choose any integral $D\in P^\circ$; under the assumption that $X$ is general in moduli, 
we know that $D$ is in the interior of the movable cone, and so is both big and movable.  Proposition 1.1 then shows 
that $X$ lies in a bounded family. \qed

\

\begin{rem}If we know that the cubic hypersurface (defined by the cubic form) in $\P ^{\rho -1} (\R )$ has connected smooth locus, 
then some part of the boundary of $\Movb (X)$ will have points on which the cubic form is strictly positive.  The above alternative argument with the volume then yields a contradiction unless there exist rigid non-movable surfaces on $X$.  Here Corollary 0.2 is true without reference to $b_3 (X)$.  Moreover, if the cubic hypersurface has two components and there are no rigid non-movable surfaces on $X$, then the same argument shows that the \Kahler cone must lie in one of the two cones corresponding to the `bounded' component of the cubic 
hypersurface and that the hyperplane $c_2 =0$ does not cut this component --- in fact the relevant cone on the bounded component is then a subcone of $\Movb (X)$.\end{rem}

The rest of this paper is devoted to applying the above theory to the case of Picard number $\rho (X) =2$, and proving Theorem 0.3.  
In particular, when $b_2 (X) =2$ and the cubic and linear forms on $H^2 (X, \Z)$ are known, 
this implies that $b_3 (X)$ is bounded.

\section{Case $\rho =2$ :  general results}

In the case $\rho =2$, the various cones under consideration, including the connected cones $P^\circ$ of the previous sections, are seen to be convex, and therefore are determined by their edge rays; this provides a significant simplification not available for higher Picard number.  Given two non-collinear 
real classes $A, B \in H^2 (X, \R)$, we shall denote by ${\rm Cone}  \langle A , B \rangle$ the closed convex cone they determine.

The proofs we give in the remaining sections were largely constructed by the author from the corresponding 
diagrams in $\R^2$, which included the real lines through the origin defined by the cubic form, any real lines 
through the origin defined by the Hessian form (a quadratic form), the line through the origin defined by $c_2 =0$ 
and rays from the origin generated by possible rigid non-movable surface classes.  We shall note below that 
appropriate real coordinates on $\R^2$ may be chosen so that things become even more concrete.
Because however such diagrams 
were needed a large number of times, it was unrealistic to include them all in the manuscript, although we do illustrate the basic three pictures (not including the line $c_2 =0$ or the possible classes of rigid non-movable surfaces) on which other relevant information may then be included for each proof.   

\smallskip
\bf To understand more easily the proofs given in the remaining sections,  the reader is strongly advised for each proof 
 to draw for themselves 
 diagrams  in $\R^2$ of the various cones, lines and rays involved.   \rm 
 \smallskip

Submerged under the details however, there is a basic strategy which the author hopes might be relevant for higher Picard number, namely using the second chern class $c_2$ to identify 
possible movable classes.  
\smallskip

We shall consider below the various basic possibilities for the cubic form on $H^2 (X, \R)$: we enumerate these here and fix on notations that we subsequently use.
 
 \
 
\noindent {\bf Case (a):}\quad  The cubic form may have three distinct real linear factors.  In this case \it real \rm (but not necessarily rational) coordinates may be chosen so that the form is $xy(x+y)$, as shown in Figure 1.  

\begin{figure}
\includegraphics[width=11cm]{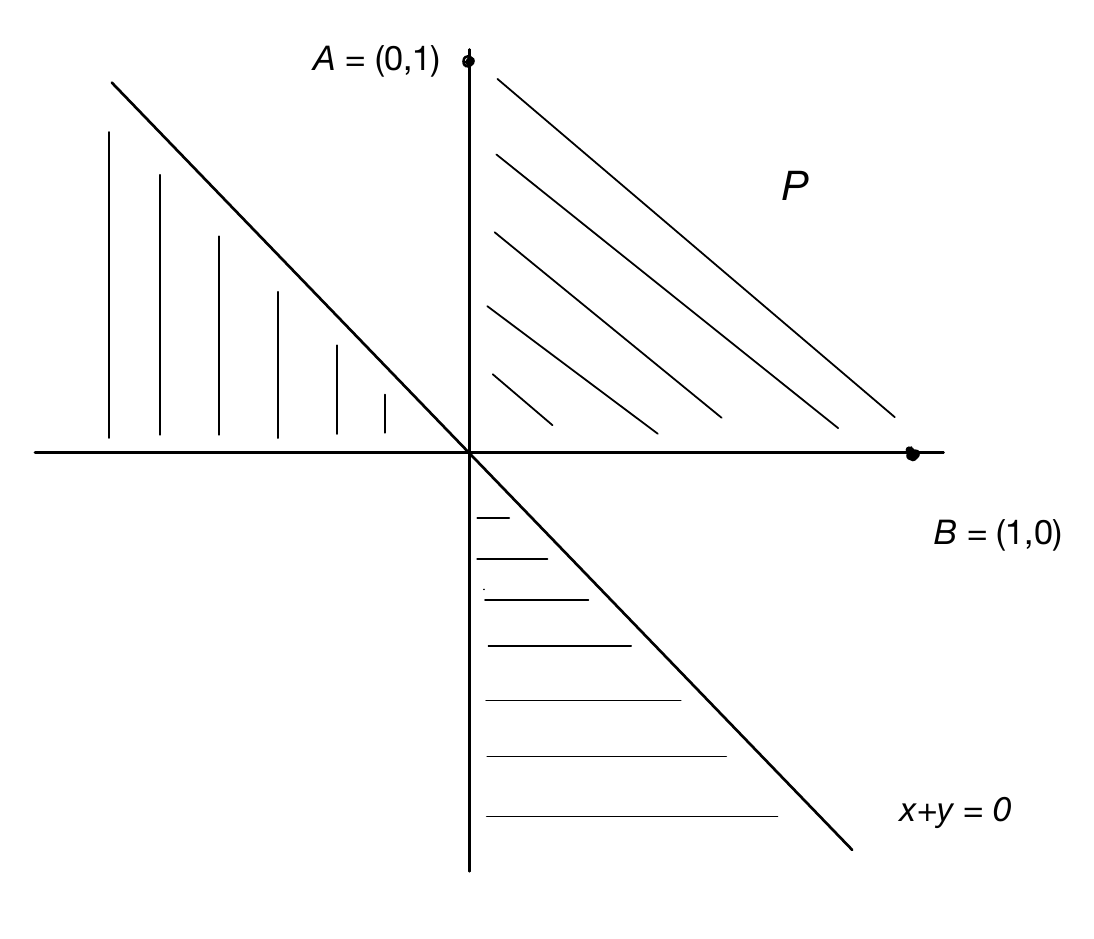}
\caption{Cubic form has three distinct real linear factors}
\end{figure}

 For any non-zero class $D= (a,b)$, the quadratic form $q_D : L \mapsto D\cdot L^2$ on $H^2 (X, \R)$, where $L =(x,y)$, has $3q_D$ given by the formula 
$(a+b)xy + ay(x+y) + bx(x+y)$, corresponding to the symmetric matrix 
$$ \left( \begin{array}{cc} b & a+b\\ a+b &  a \end{array} \right), $$ whose determinant $-a^2 - b^2 -ab$ is always negative.
Thus the index of the quadratic form $q_D$ on $H^2 (X, \R)$ 
is seen to be $(1,1)$ for all 
non-zero classes $D$ and so the only condition for $P^\circ $ is that the cubic is 
strictly positive; there are 
therefore three possibilities for $P$, the closure of the component of the positive index cone 
containing the \Kahler cone, as shown shaded  in Figure 1.

Without loss of generality we may take $P = \Cone A B$, where $A = (0,1)$ and $B= (1,0)$ --- we shall adopt this convention (when the cubic has three distinct real linear factors) for the remainder of this paper.  In the light of the two results following (Proposition 5.1 and Lemma 5.4), the class of any rigid non-movable surface $E$ lies in 
the interior of one of the two cones $\Cone {(-1,0)} {(0,1)}$ or $\Cone {(1,0)} {(0,-1)}$, since it is neither in  
$-P$ nor in the boundary of $P$.

\smallskip

\noindent {\bf Case (b):}\quad The cubic may have three (rational) linear factors  but with two being coincident, as 
shown in Figure 2. 

\begin{figure}
\includegraphics[width=10cm]{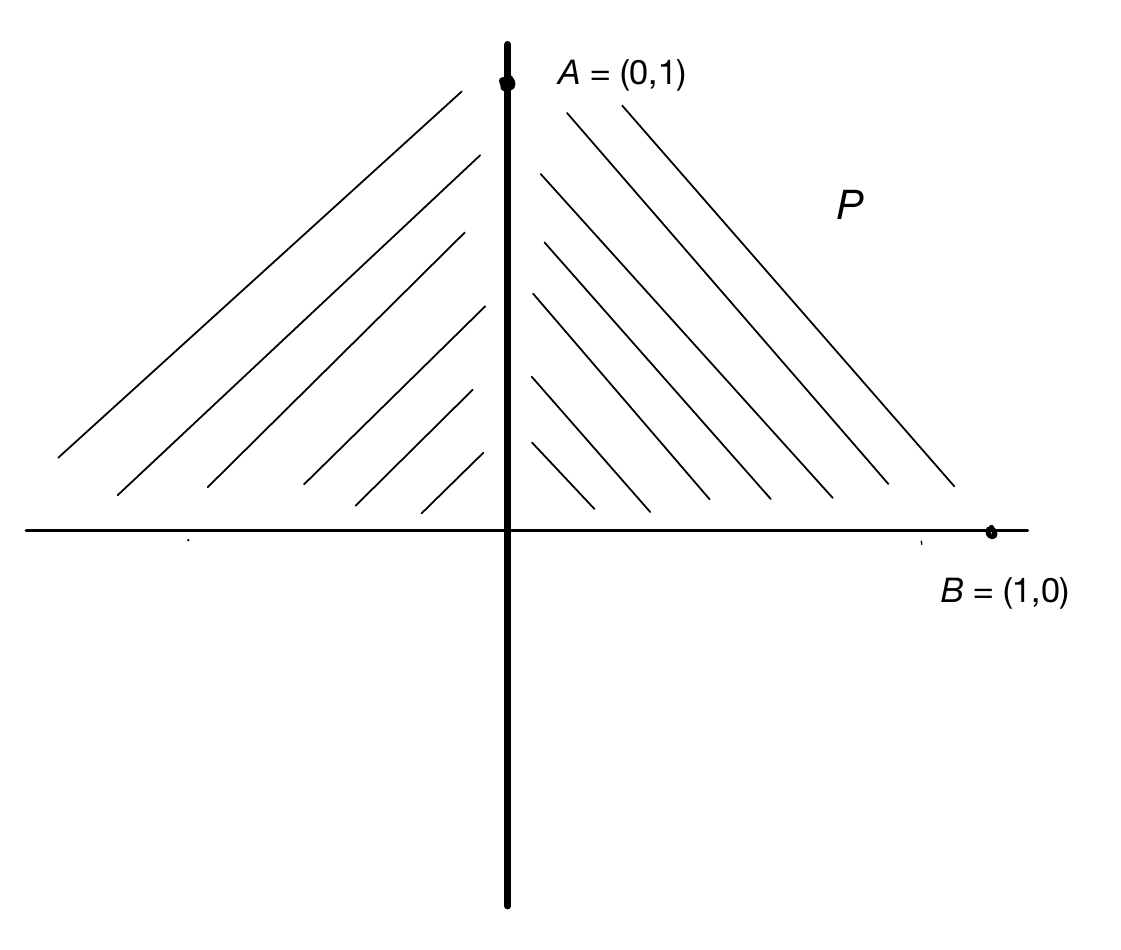}
\caption{Cubic form has a linear factor with multiplicity two}
\end{figure}

The case of three coincident linear factors cannot occur by a simple application of the Hodge Index theorem and the Hard Lefschetz theorem --- we would otherwise have linear independent  divisors $\{ D, H \}$ in $H^2 (X, \Z )$ with $H$ very ample and $D^2 \cdot H = D \cdot H^2 =0$.  From the Hodge index theorem on $H$, and the fact that $(D^2 \cdot H)H^3 = (D\cdot H^2)^2 =0$, we deduce that $D|_H \equiv 0$, a contradiction to the Hard Lefshetz theorem. Thus (rational) coordinates may be chosen so that the cubic takes the form $x^2 y$.  Here positivity is just given by $y>0$ (and $x\ne 0$) and the index condition is $x \ne 0$.  
Thus there are two possibilities for $P$, the closure of the component  of the positive index cone 
containing the \Kahler cone, which are shaded Figure 2, and without loss of generality we may take $ P = \Cone A B$ with $A= (0,1)$ and $B = (1,0)$.

Furthermore we may assume that in 
fact $A$ and $B$ are integral and primitive, at the expense of the cubic form being a rational multiple of $x^2 y$.  In this case $A$ and $B$ will generate a sublattice  of 
$H^2 (X, \Z)$ (depending only on the cubic form) of finite index.  
In the light of the two results following, the class of any rigid non-movable surface $E$ lies in 
the interior of one of the two cones $\Cone {(-1,0)} {(0,1)}$ or $\Cone {(1,0)} {(0,-1)}$.

\begin{figure}
\includegraphics[width=11cm]{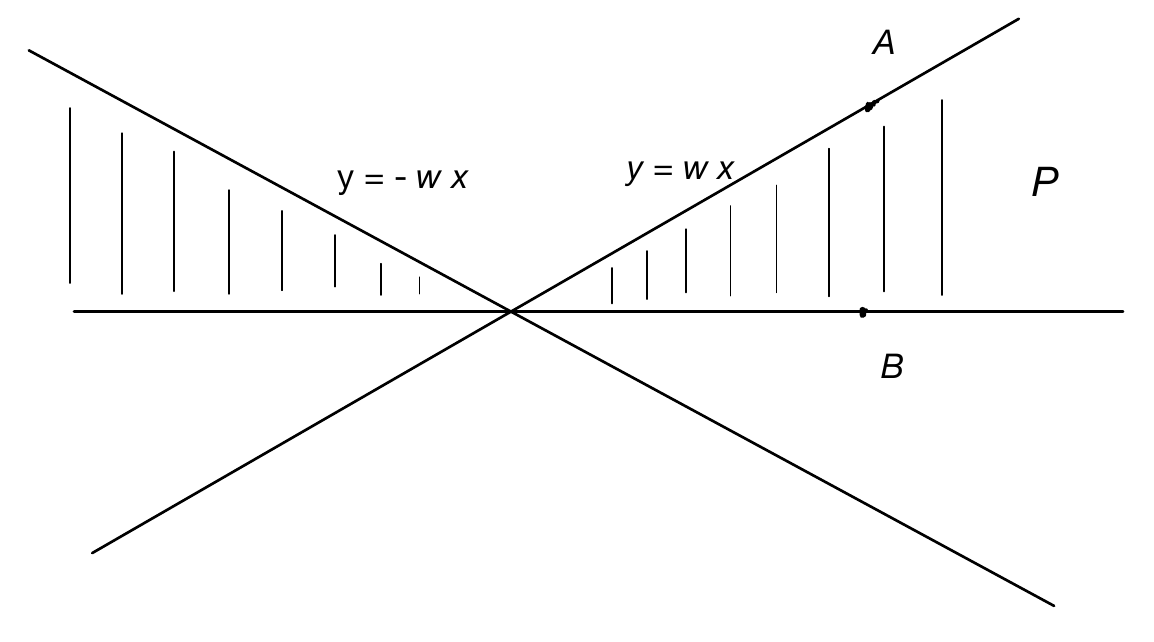}
\caption{Cubic form has only  one real linear factor}
\end{figure}

\smallskip 
\noindent {\bf Case (c):} \quad Finally we have the possibility of one real  and two complex linear factors  of the cubic.  Here real coordinates may be chosen so that the cubic is of the form $y(x^2 + y^2)$.  The positivity condition is $y>0$ and the index condition $y^2 < x^2 /3$, i.e. $|y| < |x|/\sqrt{3}$.  There are again only two possibilities for 
$P$, the closure of the component of the positive index cone  which contains the \Kahler cone, 
which are shaded in Figure 3, and without loss of generality we take $P = \Cone A B$ with $A = (1, w)$ and $B = (1,0)$, 
where $w = 1/\sqrt{3}$.  In the light of the two results following,  the class of any rigid non-movable surface $E$ lies in 
one of the regions  $\Cone {(-1,0)} {(-1, w)} \setminus \R _+ (-1,0)$ or $\Cone {(1,0)} {(1, -w)}
\setminus \R_+ (1,0)$.

\begin{prop} Suppose  $\rho (X) =2$ and $P^\circ$ is the connected component of the positive index cone which contains the K\"ahler cone; if $P$ denotes the corresponding closed convex cone,  then 
$P \subset \Effb (X)$.

\begin{proof} This follows in a similar way to the alternative proof of Theorem 0.1 given in the previous section.  We first note that it is sufficient to prove 
the result when $X$ is general in moduli, since $\Effb (X)$ at special points of moduli can 
only be larger than it is at a general point.   Suppose $P = {\rm Cone}  \langle A_1 , 
A_2  \rangle$; 
we ask about $\Movb (X) \cap P$.  This is a subcone of $P$ and therefore of the form ${\rm Cone}  \langle  B_1 , B_2  \rangle$, with $B_1 \in 
{\rm Cone}  \langle A _1, B_2 \rangle$ and $B_2 \in {\rm Cone}  \langle B_1 , A_2 \rangle$.  We claim that $A_i \in \Effb (X)$ for $i=1,2$, and hence the result follows.
By symmetry we need only prove this claim for $A_1$; if $B_1 = A_1$, then trivially $A_1 \in \Movb (X) \subset  \Effb (X)$.

Suppose now $B_1 \in P^\circ$, so in particular $B_1 \ne A_1$.  Arguing via the volume function as in the alternative proof of Theorem 0.1, we can find a nearby ray with 
integral generator $L$ for which $L \not \in \Movb(X)$ with $\vol (L) >0$, that is $L$ is big but not movable.  Therefore there exists a prime divisor $E$ (which must be rigid non-movable) such that $A_1 \in \Cone E  {B_1}$ --- in fact any $L$ as above has a $\sigma$-decomposition in the sense of \cite{Nak} as a movable class plus $bE$ for some positive real number $b$.  Finally we note from Proposition 4.4 that $E\not\in P^\circ$;  thus 
the class $A_1$ is a convex combination 
of the effective divisor class $E$ with some movable class and hence 
$A_1 \in \Effb (X)$ as claimed. \end{proof}\end{prop}

We note that we cannot have two rigid non-movable surface classes $E_i$ ($i=1,2$) such that $\Cone {E_i} A \cap P = \{ A \}$ (i.e. with $E_1 , E_2$ `on the same side' of $P$),
 since if so we have 
without loss of generality that $E_2 = aE_1 +L'$ with $L'$ ample and some $a >0$, 
 and thus that $E_2$ is big (a contradiction).
Therefore in the case $\rho =2$, there are at most two rigid non-movable surfaces on $X$. 
 In the case of two such classes, 
their classes lie with one either side of $P$ and 
no convex combination can lie 
in $-P$, since Proposition 5.1 implies that $P \subset \Effb (X)$. Moreover $\Effb (X) = \Cone {E_1} {E_2}$ in this case:  clearly $ \Cone {E_1} {E_2} \subseteq \Effb (X)$, and if it were not an equality, 
there would exist a rational class $D \in \Big (X)$ which is not in the cone 
$\Cone {E_1} {E_2}$; thus without loss of generality $E_2 \in \Cone {E_1} {D}$ and 
it is therefore big, contradicting $E_2$ being rigid non-movable.

\begin{rem} Assuming that $X$ is general in moduli, we see that if $A$ is a rational pseudo-effective class which is not in $\Movb (X)$, 
then it may be written as $L + aE$ for a well-determined rigid non-movable surface $E$, a rational class $L \in \Movb (X)$  and $a$ a positive rational number.  If we allow 
$L$ to be a real class and $a$ a real number, then there is a unique decomposition with $a$ minimal (true also when $A$ a real class), the $\sigma$-decomposition of \cite{Nak}.  \end{rem}

In the case when there are two rigid non-movable surfaces, suppose furthermore that $c_2\cdot 
E_i \le 0$ for $i = 1,2$; then both these will be zero and $c_2 \cdot H =0$ for all ample $H$ (i.e. $c_2 \equiv 0$).  This happens only when $X$ is an \'etale quotient of an abelian threefold 
\cite{Yau} and hence not 
simply connected; in any case such threefolds do form a bounded family by results of Oguiso.

In the light of Theorem 0.1 therefore, which deals with the case when there are no rigid non-movable surfaces on $X$, 
in order to prove Theorem 0.3 we may assume that $X$ is general in moduli and we 
are reduced to considering the following remaining cases.

\noindent (1)\quad There is a unique rigid non-movable surface $E$ on $X$.

\noindent (2)\quad There are precisely two rigid non-movable surfaces $E_1, E_2$ on $X$,  where $c_2 \cdot E_i \ge 0$ 
for $i=1,2$, and at least one inequality is strict.

\noindent (3)\quad There are precisely two rigid non-movable surfaces $E_1, E_2$ on $X$, where $c_2 \cdot E_1 <0$ and $c_2 \cdot E_2 >0$. 
\medskip

\begin{rem} Given the cubic and linear forms on $H^2 (X, \R)$, by Proposition 2.2 there are only finitely many possibilities for rigid non-movable classes 
$E$ with $c_2 \cdot E \le 0$, all of which have $E^3 >0$ if $c_2 \cdot E < 0$, but with a small number of possibilities with $E^3 \le 0$  if $c_2 \cdot E =0$.  
\end{rem}

In the case $\rho =2$, we can strengthen the result from Proposition 4.4.

\begin{lem}
If $E$ is the class of a rigid non-movable surface on a \CY threefold $X$ with $\rho (X) =2$,  and $P$ is the 
closure of the component of the 
positive index cone which contains the \Kahler cone, then $E\not\in P$.  

\begin{proof} For arbitrary $\rho >1$, we note by Proposition 4.4 that $E$ could only generate a ray in the boundary 
of $P$ --- in fact $E^3 =0$ 
and $c_2\cdot E =0$.  I claim that there is a mobile divisor $L$ on $X$ with $E^2\cdot L <0$.
To see this claim, we recall that $E$ corresponds to an  exceptional divisor $E'$  for a Type II or \Ttz contraction on some model $X' \ne X$ obtained by directed flops.  In the Type II case, we set $L'$ to be any very ample divisor on $X'$ and note that $(E')^2\cdot L' < 0$. In the \Ttz  case, we set $D$ to be a semi-ample free divisor on $X'$ defining the contraction and note that $(E')^2\cdot D\ < 0$; thus if $H'$ is a very ample divisor on $X'$, the very ample divisor $L' = H' + nD$ also satisfies the inequality $(E')^2\cdot L' < 0$ for $n\gg 0$.  
In both cases, we let $L$ denote the mobile big divisor on $X$ corresponding to $L'$ on $X'$.  
Note now that the formula (1) for the transformation of the cubic form under a flop yields a transformation formula for the corresponding symmetric 
trilinear form (via polarisation of the cubics); if $X'$ is achieved from $X$ after one $E$-directed flop,  this in our case reads 
$$E^2 \cdot L = (E')^2\cdot L' - ( E' \cdot\eta)^2 (L' \cdot \eta) \sum_{d>0}n_d d^3, $$ 
where the notation is as in Section 1, and hence $E^2 \cdot L <0$ as both terms on the righthand side are negative.  

Suppose now that more than one directed flop is needed to pass from the pair $(X,E)$ to $(X', E')$.
Since $E^3 =0$ 
and $c_2\cdot E =0$, it follows from the formulae ($1'$) that only two flops are in fact allowedß and there are only two possibilities how we can obtain these  numbers, 
corresponding to $9 = 8 +1 = 1+8$; both possibilities have $(E')^3 =9$ and hence $E' \cong \P^2$.
If we denote $(X', E')$ now by $(X_2 , E_2)$ and the intermediate pair by $(X_1, E_1)$, both flops must be Atiyah flops in $(-1, -1)$ curves, since the only non-zero $n_d$ 
in the formula for each flop must be one.
If we are flopping in a curve $l_2$ to pass from $(X_2 , E_2)$ to $(X_1 , E_1)$ and in a curve $l_1$ to pass from $(X_1 , E_1)$
to $(X, E)$, then by construction $E_i \cdot l_i >0$ for $i = 1,2$, and the formulae $(1')$ shows that one of these numbers is 1 and the other is 2.  
Since $-E_2|_{E_2}$ is ample, it is clear that $l_2 \not\subset E_2$, and we let $l_2' \subset E_1$ denote the flopped curve (also isomorphic to $\P ^1$) on $X_1$.  If $E_2\cdot l_2 =1$, it is just the exceptional curve of the first kind on $E_1$, and if $E_2 \cdot l_2 =2$, it is the set-theoretic double locus of $E_1$.  Since $E_1 \cdot l_2 ' <0$, it follows that $l_2' \ne l_1$.  Since $L'$ was chosen to be very ample on $X_2$, it follows that the set-theoretic base locus of the corresponding linear system $|L_1 |$ on $X_1$ is just $l_2 '$, and 
hence  $L_1 \cdot l_1 \ge 0$.
  Thus we 
may in this case just apply the above  formula twice to deduce 
that $E^2 \cdot L <0$.

We now specialize to the case $\rho =2$, and so 
$E$ generates an edge ray of $P$, and suppose $H$ is a very ample divisor on $X$.  
Since $E^2$ is not numerically trivial, it defines a non-trivial linear form which vanishes on the edge ray, and 
we deduce from Lemma 3.3 that $E^2 \cdot H >0$.  
Since $L \not \in -P$ by Proposition 5.1, it follows 
from $E^2 \cdot L < 0$ that $E$ is a convex combination of $L$ and $H$, and in particular that 
$E$ is movable, contrary to assumption.
\end{proof} \end{lem}

We finish this section with two rather technical results, which we prove now to avoid breaking the flow of later proofs.  We assume $\rho =2$ and the cubic and linear forms are known, 
along with the component $P^\circ$ of the positive index cone containing the \Kahler cone, 
with $P$ denoting its closure.  
The first of these results refers to the case where $E$ is a (known) rigid non-movable 
surface class (not necessarily unique) on $X$.  The last and  most crucial  part of this Proposition restricts how closely 
$\Movb (X)$ can approach the ray generated by  $-E$.

\begin{prop} Let $X$ be a \CY threefold with $\rho =2$, where the cubic and linear forms on $H^2 (X, \R)$ are known.
Given a rigid non-movable surface class $E$, there exists a real class (unique up to  positive multiples) $\Delta \in P$ with 
 $E\cdot \Delta ^2 =0$; moreover 
$P\cap \{ D : E\cdot D^2 \ge 0 \}  = \Cone \Delta {-E} \cap P$ and the \Kahler cone is contained in the interior of this subcone.  

If $E\cdot \Delta \equiv 0$, then $\Delta$ may be taken integral  and semi-ample, where for large enough $n$ we have a birational morphism $\phi = \phi_{n\Delta }$ contracting $E$ and with image a \CY threefold (with a canonical singularity) of Picard number one.

If $E\cdot \Delta \not \equiv 0$, then $\Delta \in P^\circ$ and 
there exists a class $R$ (with $\Delta$ and $R$ in the same open half-plane of the complement to the line generated by  $E$) such that $\Movb (X) \subset \Cone E R$, 
where the ray containing $R$ can be chosen so as to only depend 
on the cubic form, the cone $P$ and the class $E$.
 
\begin{proof} The case $E\cdot \Delta \equiv 0$ only occurs in Case (c) for the cubic, 
since then $E$ and $\Delta$ would be (rational) roots of the Hessian quadratic.

We first give a proof of the Proposition in Case (a), when the cubic form has three distinct real linear factors.  
 We adopt the notation above, so $P$ is one of the three components of the positive cone, and 
 by choice of suitable coordinates we assume that $P = \Cone A B$ with 
$A = (0,1)$ and $B = (1,0)$;  note that $E \not \in P$ by Lemma 5.4.  By Proposition 5.1, 
the class  $E$ cannot lie in the quadrant $-P$, and so it must 
lie in the interior of one of the other two quadrants; we assume without loss of generality that $E$ 
lies in the quadrant
$\Cone {(-1 , 0)} {(0, 1)}$.  

Noting that $\{ D\ :\ A^2\cdot D >0\}$ is the open right-hand half-plane and $\{D\ :\ B^2\cdot D >0\}$ is the open 
upper half-plane, our assumptions imply that $E \cdot A^2 <0$ and $E\cdot B^2 >0$; thus for $D = (x,1)$,
we know that  $E\cdot D^2$ is a quadratic in $x$ which is negative at $x=0$ and positive for $x\gg 0$, and in 
particular has a unique positive root.  Therefore 
for some $\Delta \in P^\circ$ 
(with the corresponding ray being unique) we have $E\cdot \Delta ^2 =0$.  Moreover $P\cap \{ D : E\cdot D^2 \ge 0 \} = 
\Cone \Delta B = \Cone \Delta {-E} \cap P$
and  the \Kahler cone is contained in the interior of this subcone.  
Since $\Delta ^3 >0$ and $\Delta ^2 \cdot E =0$, we then know from the index property 
at $\Delta$ that 
$\Delta \cdot E^2 <0$. We illustrate the picture in Figure 4, in which we have omitted the line $x+y =0$ so 
as not to complicate the diagram, and where the \Kahler cone is contained in the shaded subcone of $P$.  
The picture is similar for both $E^3 \ge 0$ and $E^3 \le 0$ (the case actually shown), although the proofs 
of the Proposition's last paragraph are different.

\begin{figure}
\includegraphics[width=12cm]{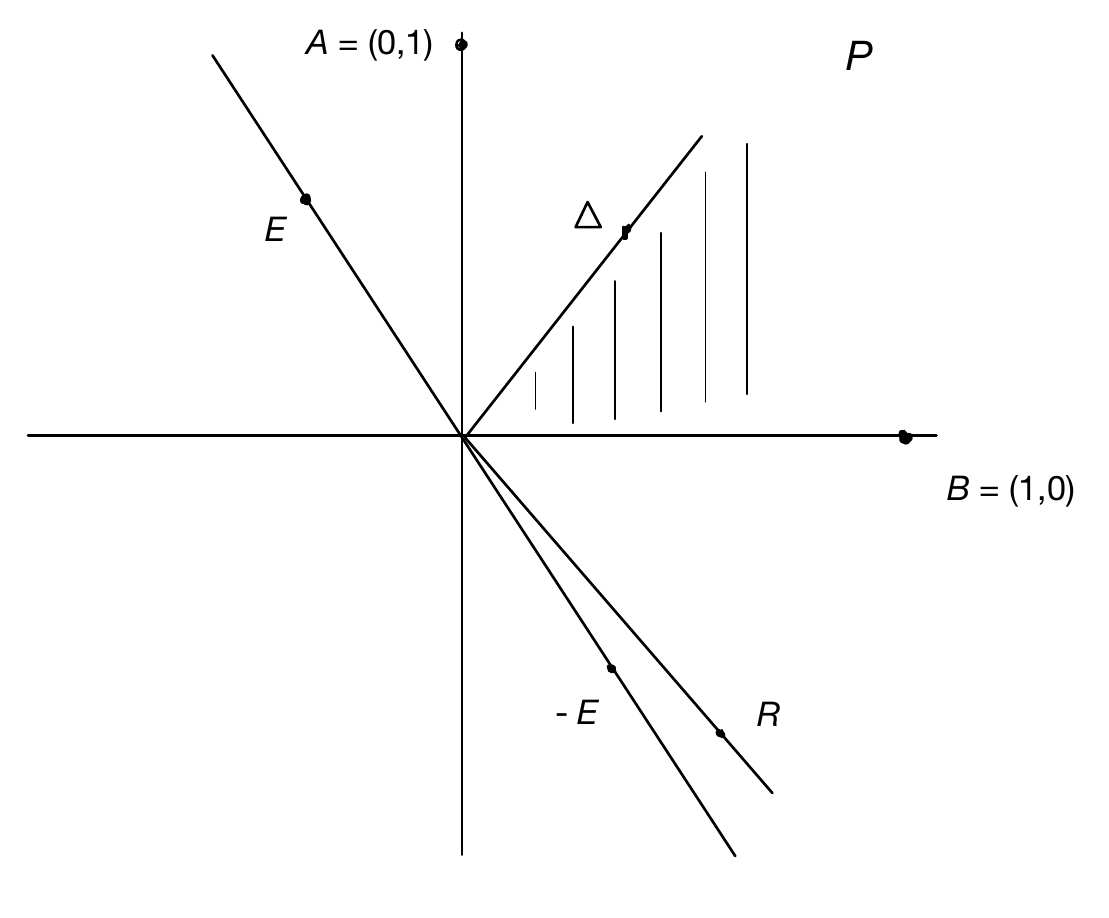}
\caption{Existence of $\Delta$ and $R$ in Case (a)}
\end{figure}

For the last paragraph, we consider separately the cases $E^3 \ge 0$ and $E^3 \le 0$.  
When $E^3 \ge 0$, we suppose that $\Delta - \lambda E$ is movable for some real $\lambda >0$ 
and bound $\lambda$ in terms of the given data.  There exists some $\mu >0$ 
(unknown) for which  the class $\Delta - \mu E$ is ample.  Therefore 
$$(\Delta - \lambda E)^2 \cdot (\Delta - \mu E) = \Delta ^3 - (\mu + 2\lambda)\Delta ^2 \cdot E + \lambda (\lambda + 2 \mu) \Delta\cdot E^2 - \lambda^2 \mu E^3 \ge 0.$$
Since $\Delta ^2 \cdot E =0$, $\Delta \cdot E^2 <0$ and $E^3 \ge 0$, we deduce that $\lambda^2$ is bounded above by $-\Delta ^3/(\Delta\cdot E^2)$.  This then gives the required statement when $E^3 \ge 0$.

When $E^3 \le 0$, we show instead that we can bound $\alpha >0$ above such that $B - \alpha E$ is movable.  In this case, we know that $B^2 \cdot E >0$, but $B\cdot E^2$ may have either sign.
There does however exist 
a $\beta > 0$  such that $B + \beta E$ is ample, and given that $E$ does not lie in $P$ and 
is a known class, 
we have a known upper bound $\beta _0$ on such $\beta$.  However for movable 
$B - \alpha E$ (with $\alpha >0$) and ample $B + \beta E$, we have 
$$ (B - \alpha E)^2 \cdot (B + \beta E) = B^3 - (2\alpha - \beta ) B^2\cdot E  +  \alpha (\alpha - 2 \beta ) B\cdot E^2  + \alpha ^2 \beta E^3 \ge 0,$$
where $B^3 = 0$ in this case.  When $\alpha > \beta /2$ the second term is negative, and by assumption $E^3 \le 0$.  If now $B\cdot E^2 \le 0$, we get a contradiction for 
$\alpha \ge 2\beta$, and so we deduce that $\alpha < 2\beta$.  If instead $B\cdot E^2 \ge 0$, then 
a similar argument shows that $B - \alpha E$ cannot be movable for $\beta/2 < \alpha < 2\beta$, and so cannot be movable for any $\alpha > \beta /2$ by convexity of the movable cone (if $D = B - 
\alpha E$ were movable for some $\alpha > 2\beta$,  we could write any element of form $B -  \alpha E$ 
with $\beta /2 < \alpha < 2 \beta$ as a convex combination of the ample divisor $B +\beta E$  and $D$). 
Thus in both cases we have 
$\alpha \le 2\beta _0 $, and this  gives the required statement when $E^3 \le 0$.

We now look at Case (b) for the cubic form, as defined above. We know that $P^\circ$ is one of the two possible components of the positive index cone, say the interior of $\Cone {(0,1)} {(1,0)}$  
(with the notation as at the start of this section).  If we have a rigid non-movable surface class $E$, then by Proposition 5.1 the class is not contained in $-P$, 
and by Lemma 5.4 it is not contained in $P$, and so in particular we have $E^3 \ne 0$.  
If $E^3 >0$, then $ E $ is 
 in the interior of $\Cone {(-1, 0)} {(0 ,1)} $, say $E= (-a, b)$ with $a>0$ and $b>0$.  Polarising the cubic form with respect to $E$ we get $x(bx-2ay)/3$.  Setting $\Delta = (2a, b) \in P^\circ$, we have $\Delta ^2 \cdot E = 0$ and the K\"ahler cone is contained in the interior of the subcone 
$P\cap \{ D : E\cdot D^2 \ge 0 \} =  \Cone \Delta {-E} \cap P$.  If $E^3 <0$, we have a similar statement, since $E = (a, -b)$ with $a>0$ and $b>0$ and $-E$ also has index $(1,1)$, so applying the 
previous calculation to $-E$ gives the required statement.  

The previous argument 
from Case (a) then yields the last sentence when $E^3 >0$.  If however $E^3 <0$, we observe that 
$A$ cannot be in  $\Mov (X)$ and we may take $R=A$  (if the rational class $A$ were in $\Mov (X)$, then we could write $A = D + H$ for 
some rational class $D \in \Mov (X)$ and some rational ample class $H$, from which we see that $A^2 \equiv 0$ is impossible).

Finally we look at Case (c) for the cubic form, where we can assume  without loss of generality 
that $P = \Cone {(1, 1/\sqrt 3 )} {(1,0)}$, as in the notation at the start of this section. 
A rigid non-movable surface $E$ cannot have class in $P$ by Lemma 5.4 and 
cannot have class in $-P$ by Proposition 5.1; in particular therefore $E^3 \ne 0$. 
If $E^3 >0$, by consideration the index of its associated quadratic form, its class  must 
either lie in the interior of the cone generated by $(-1, 0)$ and $(-1, 1/\sqrt{3})$, with the corresponding ray generated by $(-1, b)$ with $0 < b < 1/\sqrt{3}$, 
or it generates the same ray as $(-1, 1/\sqrt{3})$. 
 
Polarising the cubic form with respect to a point  $(-1, u)$ with $0 \le u \le 1/\sqrt{3}$
gives a quadratic form proportional to $3uy^2 - 2 xy + u x^2$; 
i.e. if $F = (-1, u)$ and 
$D = (x,y)$, then $F\cdot D^2$ is a positive multiple  this form.  
For $u=0$, this is just $-2xy$ (negative on $P^\circ$), and for $u = 1/\sqrt{3}$ this gives 
$(x- \sqrt{3} y)^2/\sqrt{3}$ (positive on $P^\circ$).  
For $0 < b < 1/\sqrt{3}$, we note that $3by^2 -2y +b$ is strictly positive at $y=0$, strictly negative at $y=1/\sqrt{3}$, and has a unique root inbetween.  Thus when the class of $E$ lies in the ray 
generated by $(-1, b)$ with $0 < b < 1/\sqrt{3}$,
there is a real class $\Delta \in P^\circ$ (unique up to positive multiples) with $E\cdot \Delta ^2 =0$, where again we have 
$P\cap \{ D : E\cdot D^2 \ge 0 \} =   \Cone \Delta {-E} \cap P$ and  
the K\"ahler cone is contained in the interior of this subcone.  

In the remaining case with $E^3 >0$, 
namely $E$ generating the same ray as $(-1, 1/\sqrt{3})$,
we take $\Delta$ generating the same ray as $(1, 1/\sqrt{3})$ 
and we have that $E\cdot \Delta \equiv 0$.  Note that $\Delta$ generates the wall of $P$ closest to $E$.
Therefore  $E$ and $\Delta$ represent roots of the Hessian quadratic, and so 
in particular $\Delta$ may be chosen to be integral.   
Since $\Delta = aH + bE$ for suitable ample $H$ and rationals $a, b >0$, we know that $-E|_E$ is ample.  Noting that the cubic form is strictly positive on the open upper half-plane, it follows from 
Fact 1 in \cite{WilKC} (which in turn is an easy consequence of the Key Fact in \cite{WilPic}), that
the wall of the nef cone closest to $E$ is generated by an integral semi-ample divisor $D$, defining a contraction of Types I, II or III of curves lying on $E$.  Since 
however $\Delta = \alpha D + \beta E$ for some rational $\alpha >0, \beta \ge 0$ and $\Delta |_E \equiv 0$, we deduce that $\beta = 0$, using the facts that $D\cdot C =0$ for some curve $C$ on $E$ whereas $E\cdot C = E|_E \cdot C <0$.
Thus $D$ is proportional to $\Delta$, implying that  $\Delta$ is 
semi-ample, and that for large $n$ we have a morphism $\phi = \phi_{n\Delta }$ contracting $E$ and with image a \CY threefold (with a canonical singularity) of Picard number one.

The existence of the $\Delta$ claimed in the case when $E^3 < 0$, 
where now $E \in \Cone {(1,0)} {(1, -1/\sqrt 3 )}$, follows since we may apply the previous calculation to $-E$ (which also has index $(1,1)$), where in this case $E$ cannot generate the same ray as $(1, -1/\sqrt{3})$,  since then by polarising the cubic we see that the corresponding quadratic form would be a positive multiple of $-(x-\sqrt {3} y)^2$ and hence 
strictly negative on $P^\circ$; in particular $E\cdot H^2 <0$ for any ample class $H$.

The last sentence follows in the case $E^3 >0$ with $E$ not generating the same ray as $(-1, 1/\sqrt{3})$ by
the same argument as in Cases (a) and (b), calculating with the divisor $\Delta$.  In the case $E^3 < 0$, we 
observe that $A\not\in \Mov (X)$, since there are points arbitrarily close to $A$ at which the Hessian 
is positive and therefore not in $\Mov (X)$ (cf. Lemma 3.2); therefore we may take 
$R = A$.
 \end{proof}\end{prop}
\begin{rem}
Consider  \CY threefolds $X$ with $\rho (X)=2$ and given cubic and linear forms on $\HZ 2 X$; if there is a contraction morphism $\phi _{n\Delta} : X \to Y$ of a surface $E$ to a point, then 
the Hessian form vanishes at both $E$ and $\Delta$, and  $E\cdot \Delta \equiv 0$.  Conversely given a rigid non-movable surface class $E$ on a \CY threefold $X$ which is a root of 
the Hessian form, and $\Delta$ is an integral class with $\Delta ^3 >0$
representing the other root of the Hessian, then Proposition 5.5 implies that $E\cdot \Delta \equiv 0$ 
and $\Delta$ is big and semi-ample, with the corresponding morphism contracting $X$ to a \CY 
threefold $Y$ of Picard number one,  the surface $E$ being contracted to 
a canonical singularity.
   The class of a primitive 
such $\Delta$ is determined by the Hessian form (up to ambiguity as to which 
of the two roots it corresponds).   
The fact that such threefolds (for given cubic and linear forms) constitute a  bounded family  follows
from an easy version of Proposition 1.1 (where it is unnecessary to flop).
 \end{rem}

The second technical result has the effect that for cases (2) and (3) above, knowing the class of any  integral divisor 
with some multiple being mobile is enough to yield boundedness, whilst in case (1) we may also (in the absence of bigness) 
need to specify  the class of the rigid non-movable surface.  It may be thought of as a partial generalization of Proposition 1.1 for the case $\rho =2$.

\begin{thm} We assume as above that $X$ is a \CY threefold with Picard number $\rho (X) =2$ and that the cubic and linear forms on  $\HR 2 X = \R ^2$ are given, as is the component $P^\circ$ of the positive index cone containing the \Kahler cone, with $P$ denoting its closure.
Suppose we also specify the class $L$ of a non-trivial 
integral divisor on $X$, which is not big but some positive multiple is mobile.

(i) Assuming $L$ is not semi-ample, we have boundedness for the corresponding family of \CY threefolds.
If $L$ is semi-ample, there exists at most one rigid non-movable surface on $X$.

(ii) If $L$ is semi-ample and we also know the class $E$ of a rigid non-movable surface, then we have boundedness of the family.

\begin{proof} (i) Recall that $L \not \in -P$ by Proposition 5.1.  Since some multiple $mL$ is mobile, we may apply Bertini's theorem to deduce that $mL$ is linearly equivalent to either 
an irreducible divisor or a positive multiple of one.  In both cases, $L$ corresponds on some minimal model $X'$ to a nef divisor $L'$, which is semi-ample by \cite{Og}. 
 As the movable divisor $L$ is not big,  we must have $(L')^3 =0$, and hence by repeated use of formula (1) we have 
 $L^3 \le 0$ on $X$.
 If $L \not \in P$, we can choose any integral divisor $D'$ in the interior of 
 the cone generated by $L$ and $P$ but not in $P$, 
 which will then be movable and big (for any ample $H$, $ \Cone H L \subset \Movb (X)$), and so boundedness follows from Proposition 1.1.

We are reduced therefore to the case when $L$ generates an edge ray of $P$ and that $L^3 =0$.  By the transformation formulae (1) 
 for the cubic form under flops and the assumption that $L$ is not big, we deduce that $L$ is already nef and  hence semi-ample and is in a wall of the nef cone.  
   
 There are then two cases to consider, namely  $L^2 \not \equiv 0$ and $L^2 \equiv 0$, corresponding to the 
 morphism $\phi _{nL}$ for $n$ sufficiently large defining an elliptic fibre space structure on $X$, respectively a K3 or abelian fibration over $\P ^1$.  We consider the half-plane  (containing $P$) with $L$ in its boundary.
Suppose now that $E$ is a non-movable rigid surface on $X$.  Since $E \not\in -P$ by 
Theorem 5.1, we note that $E$ is in this half-plane (otherwise, we could write $L$ as a convex combination of $E$ and an 
ample divisor $H$, contradicting the assumption that $L$ is not big).  As $E$ cannot be a multiple of $L$, it lies in the open half-plane.  If $L^2 \not \equiv 0$, the half-plane is defined by the linear form corresponding to $L^2$ and so $L^2 \cdot E >0$.  If $L^2 \equiv 0$  and $H$ is any ample divisor, the half-plane is defined by the linear form $L\cdot H$ and so $L\cdot E \cdot H >0$.

If there was another rigid non-movable surface class $E'$, it would have to lie 
`on the other side' of $P$ to $E$, i.e. in the complementary half-plane; this is not possible 
by the argument just given.
Thus $E$ is the only rigid non-movable surface class and is `on the other side' of $P$ to $L$, in other words  $nL + E$ 
 is ample on $X$ for all $n\gg 0$.  We note for later use that 
 in both cases 
 $$h^0 (X, \O _X(nL -E)) =0\ \hbox{\rm for all\ }  n>0 $$ (just intersect with $L^2$, respectively 
 $L \cdot H$).

(ii)  We assume that the class $E$ of the unique rigid non-movable surface is also given.  We 
deal first with the case where there is an elliptic fibre space  
 structure on $X$ 
 defined by a morphism $\phi _{nL} : X \to W$ for $n$ sufficiently large.  Moreover $L^2 \cdot E >0$ is known.  We show that this is enough information to give not just birational boundedness (which is implied by \cite{GrEll}) but also actual boundedness.
 
 We know that the base surface $W$ is $\Q$-factorial with only finite quotient singularities, and is rational as $X$ is assumed to be simply connected (\cite{OgFin}, \S 3). Moreover, results 
 of Alexeev \cite{Al} imply that these base surfaces form a bounded family (\cite{OgFin}, \S 5 or \cite{OPpos}, Proof of Theorem 1).  Since in our case $\rho (W) =1$, there is a universal positive integer $r$ such that for any effective Weil divisor $H$ on one of our base surfaces $W$, we have that $rH $ is Cartier and very ample.   In particular we deduce on $X$ that $ L' = rL$ 
 is free, for some $r$ independent of $X$ and the elliptic space structure, albeit $r$ not explicitly given.  We also have $h^0 (X, \O _X(nL' -E)) =0$ for all $n>0$, and so
 we deduce that $h^0 (X, \O _X(nL' )) \le h^0 (E, \O _E(nL'))$.  Since $L'|_E$ is free, we can argue via the short exact sequence 
 $$ 0 \to \O_E (((n-1)L') \to \O _E (nL') \to \O_C (nL') \to 0,$$
 (where $C = L'|_E$ is a smooth irreducible non-rational curve) and the estimate $h^0 (C, \O_C (nL')) \le  nL'\cdot C = n r^2 L^2 \cdot E$, to deduce by induction an effective upper-bound 
 (namely a quadratic formula in $n$ whose coefficients depend on the universal integer $r$ as well the given invariants, 
 the class $E$ and the cone $P$)
 for $h^0 (E, \O_E (nL'))$ and hence for $h^0 (X, \O_X (nL'))$. 
 
 Since however we also know $E$, there exists a fixed $s>0$ for which $D= sL' +E \in P^\circ$ satisfies $D^2 \cdot E >0$.  Assuming without loss of generality that $X$ is general in moduli, 
 it follows from Lemma 4.3 and Proposition 4.1 that there exists 
  $m>0$ depending only on $D^3$ and $c_2 (X)\cdot D$  with 
 $h^0(X, \O _X(mD))>1$.  The same argument using the Riemann--Roch theorem also implies that 
 $$h^0(X, \O _X(nmD)) \ge \chi (\O _X(nmD)) = {1\over 6} n^3m^3 D^3 + {1\over 12} nm D \cdot c_2 (X),$$
 and hence we have  
 a cubic formula  in $n$ (coefficients depending 
 the universal integer $r$ as well as the given invariants,  the classes $L$ and $E$ and the cone $P$) giving a lower bound for 
 $h^0(X, \O _X (nmD))$. Thus for some $n$, also depending only on the universal integer $r$ as well as the given invariants,  the classes $L$ and $E$ and the cone $P$, we have $h^0(X, \O _X (nmD))
 > h^0(X, \O _X (nmsL'))$.
 Therefore the mobile part $M$ of $|nmD|$ is $|nmsL' + cE|$ for some $0 < c \le nm$.  
 We note that $NL + E$ is ample for all $N\gg 0$, and so there exists an ample divisor class  
 $H \in \Cone {L}{M}$; therefore $M$ is a convex combination of $H$ and $E$, and hence
 $M$ is also big.
 
 Since $c \le nm$, there are only finitely many possibilities for the big mobile divisor class $M$ 
 found in this way, where the bounds on $M^3$ and $M\cdot c_2$ depend not only on the initial data but also on the fixed but non-explicit integer $r$. Boundedness follows from Proposition 1.1.

The second case is where the morphism $\phi _{nL}$ for $n \gg 0$ is 
 a K3 or abelian fibre space over $\P ^1$.  Here we note that $L^2 \equiv 0$, and so $(xE + zL)^3 = x^2 y$, where $y = (E^3 x + 3L\cdot E^2 z)$;  
 the cubic must then in appropriate rational coordinates be as in Case (b) above, with $L$ in the edge ray of $P$ generated by $A$.
 Since $NL+E$ is ample for all $N\gg 0$, we note that $L\cdot E^2  = L\cdot (NL +E)^2 >0$; suppose $s>0$ has the property  that $sL +E$ is ample (in the following argument we do not need to 
 know what $s$ is), and then $H = 10(sL +E)$ is very ample by \cite{Og}.  For any $n>0$, $h^0 (X, \O_X (nL-H)) =0$, since $L\cdot (nL-H)\cdot H = -L\cdot H^2 <0$.
  Taking cohomology of a standard short exact of sheaves then yields $h^0 (X, \O_X (nL)) \le h^0 (H, \O_H (nL))$ for $n>0$.  Now apply the same argument on $H$, 
 observing that $h^0 (H, \O_H (nL -H)) =0$ since $(nL-H)\cdot H \cdot (NL +E) \sim -100 N L\cdot E^2 <0$ for $N$ suitably large compared with $n$ and $s$.  Therefore  
 $h^0 (H, \O_H (nL)) \le h^0 (C, \O_C (nL)) \le n L \cdot H^2 = 100n L\cdot (sL+E)^2 =
 100 n L\cdot E^2$ for $n>0$, where the smooth curve $C$ denotes the intersection of two generic elements of $|H|$.  Thus 
 $h^0 (X, \O_X (nL)) \le 100 n L\cdot  E^2$.    
 
 However, given the classes of $L$ and $E$,  we can find a fixed $r>0$ for which $D= rL +E \in P^\circ$ satisfies $D^2 \cdot E >0$.  Assuming again that $X$ is general in moduli, it 
 follows from Lemma 4.3 and Proposition 4.1 that we can find a positive integer $m>0$ depending only on known data  with 
 $h^0(X, \O _X(mD))>1$.  As in the previous case, the argument using the 
 Riemann--Roch theorem yields  
 a formula (cubic in $n$ with coefficients depending only on known data) giving a lower bound for 
 $h^0(X, \O _X (nmD))$ for $n>0$. 
 Therefore we can find a positive integer  $n$ depending only on the known data 
 (and not on the integer $s$) with 
 $h^0 (X , \O _X (nmD)) > h^0 (X , \O _X (nmrL))$.
 Hence the mobile part $M$ of $|nmD|$ is $| nmrL +cE|$ for some $0<c \le nm$, and it is  big for the reasons given in the previous case.   Arguing  as in 
 the previous case, there are only finitely many possibilities for the big mobile divisor class $M$ 
 found in this way, 
 and we may again apply Proposition 1.1 to deduce boundedness.
\end{proof} \end{thm}

\section{Proof of Theorem 0.3 for cases (1) and (2)}

Recall from the previous section that case (1) was when there was only one rigid non-movable surface on $X$,  and 
case (2) was when there were two, $E_1$ and $E_2$ with $c_2 \cdot E_i \ge 0$ for $i=1,2$, with at least one inequality strict.  Throughout we can assume that $X$ is general in moduli.

\begin{prop} Consider the set of \CY threefolds $X$ with Picard number $\rho (X) =2$ 
and fixed cubic and linear forms on $\HZ 2 X$, for which 
either case (1) or case (2) applies.  Suppose we are given  finitely many integral classes $M_1, \ldots , M_N 
\in \HZ 2 X$ (depending on the given data) 
such that if $X$ is a \CY threefold in the above set, then for at least one $i$ we have 
$h^0 (X, \O_X (M_i)) >1$ (with $M_i$ also big in case (1)); then the set of such $X$ forms a bounded family.

\begin{proof} 
Suppose we have $X$ a \CY threefold as above;  we can assume that  for $M= M_i$ say,   we have $h^0 (X, \O_X (M)) >1$ (with $M$ big for case (1)).

We first prove the Proposition in case (1), 
 where the unique rigid rational surface class (as yet unknown) will be 
denoted by $E$.  Assuming $X$ is general in moduli, we can 
write $|M| = |L| + aE$ where some multiple of $L$ is mobile and $a\ge 0$ an integer (cf. Remark 4.2).  
Recall that $c_2 \cdot L \ge 0$ by Remark 1.2.
If $c_2 \cdot E \le 0$, then by Proposition 2.2 there are only finitely many possibilities for the class $E$, and for each possibility, using the last part of Proposition 5.5 we see that 
the possible values of $a$ have an upper bound 
(depending only on the given data and the choices for $M$ 
and $E$).  If however $c_2 \cdot E >0$, then either $a=0$, or both $a$ and $c_2 \cdot E$ are bounded, the latter fact implying by Proposition 2.2 that there are only finitely many possibilities for $E$.   Thus either $M$ is movable and big, in which case boundedness follows from Proposition 1.1, or there is a finite set of possible classes for $E$ and the class $L$, and hence by Proposition 1.1 and Theorem 5.7 the resulting family of threefolds $X$ is bounded.

When case (2) applies, we denote the two rigid non-movable surfaces by $E_1$ and $E_2$, where we 
assume that $c_2 \cdot E_1 >0$ and $c_2 \cdot E_2 \ge 0$.  We write $|M| = |L| + aE_1 + bE_2$, where some multiple of $L$ is mobile; 
since $c_2 \cdot L \ge 0$, we deduce that either $a=0$, or that both $a$ and $c_2 \cdot E_1$ are bounded, the latter implying that there are only finitely many possibilities for $E_1$
by Proposition 2.2.  Thus there are only finitely many possibilities for $aE_1$.  If $c_2 \cdot E_2 >0$, then the same argument says that there are only finitely many possibilities for $bE_2$.  If however $c_2 \cdot E_2 =0$,  there are only finitely many possibilities for $E_2$
by Proposition 2.2, and then $b$ is bounded 
by the last part of Proposition 5.5.  Thus there are only finitely many possibilities for the movable class $L$, and so the set of threefolds is bounded by Proposition 1.1 and Theorem 5.7 (i) (where in this case $L$ cannot be semi-ample as there are two rigid non-movable surface classes).
\end{proof} \end{prop}

\begin{thm} If $\rho (X) =2$, the cubic and linear forms are known and there is exactly one 
(unknown) rigid non-movable surface $E$, then boundedness holds.

\begin{proof} Throughout the proof, we may assume that $X$ is general in moduli.  
We may also assume that the component $P^\circ$ of the positive index cone that contains the 
\Kahler cone has been specified, with closure denoted by $P$.
Given these data, we prove that finitely many integral classes $M$ may be found, where for each $X$ 
at least one such $M$ will be big with $h^0 (X, \O_X (M)) >1$.  
The result follows by applying Proposition 6.1.  

To exhibit these divisors $M$,  we shall find a finite set of integral divisors $D \in P^\circ$ such that,  
for a given threefold $X$, one of these classes will satisfy  $D^2\cdot E \ge 0$.  Then by Lemma 
4.3, for a very ample divisor $H$ generic in its linear system, we have $D|_H$ is nef.  Applying 
Proposition 4.1, $D$ is big and 
we can find an integer $m>0$, depending on the given data and $D$, such that 
$h^0 (X, \O_X (mD)) >1$.  For such a divisor $D$ we can set $M = mD$.  
In the proof below, there will be a finite number of cases, and in each case there may be a finite number of subcases, for instance 
involving the choice  of the class $E$ from a finite set of possibilities, and 
in each of these subcases we will show how to choose a suitable divisor class $D$ which works.

(a)  Consider first Case (a) for the cubic, as defined in the previous section; we may assume that coordinates have been chosen such 
that $P = \Cone {A} {B}$, where $A= (0,1)$ and $B= (1,0)$.
As in the proof of Proposition 5.5, we may without loss of generality take $E$ to be in the interior of 
$\Cone {(-1,0)} {(0,1)}$, switching the two coordinates if necessary. 
We consider the two subcases, namely $c_2 \cdot E \le 0$ and $c_2 \cdot E >0$.  

In the former subcase, the class $E$ is determined up to finitely many possibilities
by Proposition 2.2, so we may assume it is known. We apply Proposition 5.5 yielding a 
divisor class $\Delta \in P^\circ$, unique up to positive multiples, for which $\Delta^2\cdot E =0$.  
 Now choose an integral $D$ in the interior of $\Cone \Delta B$; having made the choice, we 
may also regard $D$ as part of the known data and note that $D^2\cdot E >0$, so that $D$ is the divisor we seek.

The other subcase that we need to consider for the proof of the theorem for (a) is when we have $c_2 \cdot E >0$.  We may assume that $E^3 < 0$, 
since otherwise by Proposition 2.2 there are again only finitely many possibilities for $E$ and the previous 
argument applies.  Since $E\not\in P$ by Lemma 5.4,  $E$ may be taken 
in the interior of $\Cone {(-1 ,1)} {(0,1)}$ and $\Delta$ is then in the interior of  $\Cone {(0, 1)} {(1,1)}$, and so for any choice of integral $D$ in the interior 
of $\Cone {(1,1)} {(1,0)}$, we have $D^2 \cdot E >0$ and so $D$ is the divisor we seek.

We now essentially repeat this argument for Cases (b) and (c)
for the cubic form (as described in Section 5), indicating where changes are needed.  

(b) In Case (b) for the cubic we know that $P^\circ$ is one of the two possible components of the positive index cone, which we may take to be the interior of $\Cone {(0,1)}{(1,0)}$.  If we have a rigid non-movable surface class $E$, then  
 $E$ is not in $-P$ by Proposition 5.1,  and it is not in $P$ by Lemma 5.4, so $E^3 \ne 0$.  
  If $E^3 > 0$, then $ E $ is in the interior of $\Cone {(-1 , 0)} {(0 ,1)} $, with $E= (-a, b)$ where $a > 0$ and $b>0$.     By Proposition 2.2 there are only finitely many possibilities for the class $E$ such 
  that $E^3 > 0$, and for each of these we have $\Delta \in P^\circ$ with $\Delta ^2 \cdot E =0$.
  Again therefore we may just choose an integral $D \in \Cone {\Delta} B$, which therefore satisfies 
  $D^2 \cdot E >0$, and this is the divisor being sought.
  
We assume therefore that $E^3 <0$, and so $E = (a, -b)$ with $a>0$ and $b>0$.
We noted in the proof of Proposition 5.5 that $A \not\in \Mov (X)$
(since otherwise we could write $A = D + H$ for 
some rational class $D \in \Mov (X)$ and some rational ample class $H$, from which we see that $A^2 \equiv 0$ is impossible),
 but $A \in \Effb (X)$ by Proposition 5.1.
Since $E$ is the unique rigid non-movable surface on $X$, any effective class is the sum of a movable class and a non-negative multiple of $E$.  
 We deduce by Proposition 5.1 that $A$ in this case 
determines a boundary ray of both the movable and pseudo-effective cones.  In particular, 
from the former statement and Remark 1.2  we 
have $c_2 \cdot A \ge 0$.

If $c_2 \cdot E \le 0$, then by Proposition 2.2 there are finitely many possibilities for $E$, 
and for each possibility we have $\Delta \in P^\circ$ with $\Delta ^2 \cdot E =0$.
We choose  an integral divisor $D$ in the interior of $\Cone A \Delta$, and hence 
with $D^2\cdot E >0$; and this is the divisor being sought.
 
We now need to cover the case $E^3 <0$ and $c_2 \cdot E >0$.   We noted above that  
$c_2 \cdot A \ge 0$, and so 
the line $c_2 =0$ 
cannot cut $P^\circ$,  since otherwise $c_2 \cdot E \le 0$.  Similarly the line 
$c_2 =0$  cannot contain $B$.  
We suppose first that the line \it does not \rm contain $A$.    As the 
line $c_2 =0$ does not now cut $P$, we may take integral $F$ in the interior of $\Cone {(0, -1)} {(1, 0)}$ generating the line
$c_2 =0$, and the argument from Proposition 5.5 produces a real class 
$\Delta \in P^\circ$, unique up to positive multiples, 
such that $\Delta ^2 \cdot F =0$, i.e. $c_2$ is a positive multiple of $\Delta ^2$.  Thus taking an integral $D$ in the interior of $\Cone A {\Delta}$, we note that $D^2\cdot E >0$ 
for all possibilities for the rigid non-movable class $E$, and so $D$ is the divisor we seek.

The only difficult case is when $E^3 <0$ and $c_2 \cdot E >0$,  and $c_2 =0$ \it does \rm contain the primitive  
class $A$; here we exploit the integral structure to find 
an integral class $D \in P^\circ$ for which $D^2 \cdot E >0$ 
for all the possible rigid non-movable surface classes $E$ with $c_2 \cdot E >0$.  Recall that the cubic form is written as $k_1 x^2 y$ for some positive rational $k_1$ and that $A$ and $B$ generate a sublattice  of 
$H^2 (X, \Z)$ (depending only on the cubic form) of finite index; 
the linear form $c_2$ 
is now of the form $k_2 x$ for some positive rational $ k_2$.  
Suppose that $E$ is a rigid non-movable surface on $X$, then by Lemma 2.4 we have $E^3 + 
{1\over 2} (c_2(X)\cdot E)^3 \ge - 99$.
Thus if 
we write $E = (a, -b)$, where $a,b$ 
will be positive rationals whose denominators are bounded above, this says that 
$-k_1 a^2 b +  {1\over 2}(k_2 )^3 a^3  \ge - 99$. If we set $c = (k_2)^3/k_1$, we have that 
$$b \le {1\over 2} ca + 99/(k_1 a^2).  \eqno{(4)}$$
Given that the denominator of $a$ is bounded, there are only a finite number of possible values of $a$ for which 
$ k_2 ^3 a^3 \le 198$, i.e. such that $99/(k_1 a^2) \ge {1\over 2} ca$.  Thus there are only a finite number of possible values of $a$ for which ${1\over 2} ca + 99/(k_1 a^2) \ge ca$.  Since the denominator of $b$ is also bounded, for each value of $a$ there will only be a finite number of 
possible values of $b$ for which the inequality (4) holds.  It follows that 
there will
only be finitely many possible classes $(a, -b)$ of rigid non-movable surfaces for which $b > ca$, and thus for some appropriate $c' \ge c$, depending only on the known quantities, we will have $b \le c' a$ for all 
possible rigid non-movable surface classes.  If $F \in \Cone B {-A}$ lies on the line $y+ c' x =0$, then we take $\Delta \in P^\circ$ with $\Delta ^2 \cdot F =0$ and choose an integral 
$D $ in the interior of $\Cone A \Delta$.
Then $D^2 \cdot E > 0$ for all the possible rigid non-movable classes $E$ in this case, and so $D$ is the divisor we seek.

(c)  In Case (c) for the cubic, as defined in the previous section, we may assume that $P = \Cone A B$ with $A= (1, 1/\sqrt 3 )$ and $B= (1,0)$
as in previous conventions.  If $E^3 \ge 0$, then its cohomology class 
must be in  $\Cone {(-1, 0)} {(-1, 1/\sqrt 3)}$ by consideration of the index, and 
by Proposition 5.1 it is not in the ray generated by $(-1,0)$.  By Remark 5.6, we may also assume that $E$ is not in the ray generated  by $(-1, 1/\sqrt 3 )$,
since in that case we already have boundedness.  Thus $E$ is in the 
interior of the cone $\Cone {(-1, 0)} {(-1, 1/\sqrt 3)}$.  Since now
$E^3 >0$,  
by Proposition 2.2 there are only finitely many possibilities for the class of such an $E$, and for each possibility there exists 
a  real class $\Delta \in P^\circ$ with $\Delta ^2 \cdot E =0$, and we then 
 just choose an integral class $D$ in the interior of $ \Cone {\Delta} B$.  As before $D^2 \cdot E >0$, 
 and so $D$ is the divisor we seek.

Finally, we show that $E^3 <0$ does not occur in this case; if it did, its class would be in 
the interior of the cone 
$\Cone {(1,0)} {(1, -1/\sqrt 3 )}$ --- if $E$ were a multiple of $(1, -1/\sqrt 3 )$, then 
for $Z = (x,y)$, we see by polarising the cubic that 
$E\cdot Z^2$ is a positive multiple of $- (x- \sqrt 3 y)^2$, which 
would be negative on $P^\circ$, and hence we would have $E\cdot H^2 <0$ for any ample divisor $H$. 
Now $A\not \in \Mov (X)$,
since there are points arbitrarily close to $A$ at which the Hessian 
is positive and therefore not in $\Mov (X)$ (cf. Lemma 3.2), 
 but $A \in \Effb (X)$ by Proposition 5.1.  
Since $E$ is the unique rigid non-movable surface on $X$, any effective class is the sum of a movable class and a non-negative multiple of $E$.
We deduce that $A$ 
in this case determines a boundary ray of both the movable and pseudo-effective cones
 (cf. the proof for 
part (b)).  As in the alternative proof of Theorem 0.1
from Section 4, we note that 
 $\vol (D) \ge D^3$ for rational classes $D\in \Mov (X)$ and so by continuity that  $\vol (A) \ge A^3$. 
Again using continuity of $\vol$, we can find a nearby ray 
 with 
integral generator $L$ for which $\vol (L) >0$ but $L \not \in \Effb(X)$; such an $L$ would be big, providing the required contradiction.\end{proof}\end{thm}

If there are two rigid non-movable surfaces $E_i$ ($i=1,2$) on $X$, then we are in cases (2) or (3) from Section 5.  Throughout we shall assume that $X$ is general in moduli.  Moreover we showed just after Proposition 5.1 that $\Effb (X) = \Cone {E_1} {E_2}$. The rough idea behind the proof of the next theorem is as follows: under the assumptions, we can reduce to the case where a 
 divisor $D_2 \in P^\circ$ may be found, unique up to positive multiples, with $c_2$ a positive multiple 
of $D_2^2$.  If the ray generated by $D_2$ is rational, we may take $D_2$ to be integral, and then it is this divisor we use to prove boundedness.  If the ray  is not rational, then the argument is less easy, 
and we have to deal with various other subcases before being able to employ the main argument 
using fractional divisors which will complete the proof.  This latter argument will only be needed when 
the cubic form has three distinct real factors (Case (a)).
\medskip

\begin{thm}  The result in Theorem 0.3 is true in case (2) of Section 5.
\begin{proof}
 By assumption we have rigid non-movable surfaces $E_i$ with $i = 1,2$, with  $c_2 \cdot E_i  \ge 0$ and the inequality is strict 
for at least one of them.
Throughout the proof, we may assume that $X$ is general in moduli.  
We may also assume that the component $P^\circ$ of the positive index cone that contains the 
\Kahler cone has been specified, with closure denoted by $P$.
Given these data, we prove that finitely many integral classes $M$ may be found, where for each 
case (2) threefold $X$ with the given data, 
at least one such $M$ will have $h^0 (X, \O_X (M)) >1$.  
The result then follows by applying Proposition 6.1. 

In most of the subcases listed below, we shall find the divisors $M$ by 
exhibiting finitely many integral divisors $D \in P^\circ$ such that,  
for a given threefold $X$, one of these classes will satisfy  $D^2\cdot E_i \ge 0$ for $i=1,2$.
  If this is true for a given divisor $D$ on some $X$ with the given data, then by Lemma 
4.3, for a very ample divisor $H$ on $X$,  generic in its linear system, we have $D|_H$ is nef 
and big.  Applying 
Proposition 4.1, 
we find an integer $m>0$, depending on the given data and $D$, such that 
$h^0 (X, \O_X (mD)) >1$.  For such a divisor $D$ we may set $M = mD$.

(a) We consider first the possibility Case (a)  for the cubic (where there are three distinct linear factors), where as before we may take 
$P = \Cone A B$, with $A= (0,1)$ and $B= (1,0)$.  By Proposition 5.1 and Lemma 5.4, one of the $E_i$ must lie in the 
interior of the quadrant $\Cone {-B} A$ and the other in the interior of the quadrant 
$\Cone B {-A}$.
Our assumptions then imply that the linear form $c_2$ does not vanish on a wall of $P$ and does not cut the interior of  $P$; it is therefore either the third line of the cubic (that is 
given by $x+y =0$ in the coordinates chosen before) or without loss of generality it cuts the interior of the cone $\Cone B {(1,-1)}$.  We let $F_1$ denote a  positive 
multiple of $(-1,1)$ and 
$F_2$ the primitive integral generator 
in $\Cone B {(1,-1)}$ of the line $c_2 =0$; in the case when $c_2 =0$ is the third line of the cubic, we shall take $F_1 = -F_2$.  Exactly as in the proof of the first paragraph of Proposition 5.5 (but with $F_i$ rather than $E$), we can choose  
 corresponding classes 
$D_i$  in the interior of $P$, unique up to positive multiples, for which $D_i ^2 \cdot F_i =0$ ($i=1,2$); for instance $D_1$ lies in the ray generated by $(1,1)$.   In particular note that  $D_2 \in \Cone {D_1} B$ and that $c_2$ is a positive multiple of $D_2^2$.

In the exceptional case where $c_2 =0$ is the third line of the cubic, we 
may take $D_1 = D_2$ and denote this class by $D$.  Here the equation $F_2 \cdot Z^2 =0$ defines two lines, one of them generated by the integral class $F_2$ and one generated by 
$D \in P^\circ$.  
Both lines are therefore defined over the rationals, and we may take $D$ to be integral in this case.  Note that the conditions on the $c_2 \cdot E_i$ imply that $D^2 \cdot 
E_i \ge 0$ for $i = 1,2$, with strict inequality for at least one.  
Thus $D$ is an appropriate divisor in this subcase.

In the general case (where $c_2 =0$ is not the third line of the cubic), we may assume that $E_2 \in \Cone B {F_2}$ (no longer assuming that $c_2 \cdot E_2$ is  strictly positive).  By Lemma 5.4,  it is not in the ray generated by $B$, and in particular has $E_2 ^3 <0$; 
on the other hand $E_1 \in \Cone {-F_2} A$ may have either $E_1 ^3 \ge 0$ 
or $E_1 ^3 <0$.  We let $\Delta _1 , \Delta _2$  be the (as yet unknown) classes (unique up to positive multiples) in 
the open quadrant $P^\circ$ determined by $E_1 , E_2$ as in the first part of Proposition 5.5.
The easy subcase is when $E_1 ^3 < 0$, since then $\Cone {D_1} {D_2} \subset \Cone {\Delta _1} {\Delta _2}$.
We may just choose an integral divisor $D$ in the interior of $\Cone {D _1} {D _2}$, which ensures that $D^2 \cdot E_i >0$ 
for $i = 1,2$. Thus $D$ is the divisor we seek in this subcase.

The more difficult subcase is when  $E_1 ^3 \ge 0$; here  
there are only finitely many possibilities for its class, which we may now assume 
 is known. The class $\Delta _1$ is then also known, up to positive multiples, and 
 $P^\circ \cap \{ E_1 \cdot L^2 >0\}$ is the interior of 
 $\Cone {\Delta_1} B$.  Suppose now that $c_2 \cdot E_1 >0$; since $D_2^2$ is a positive multiple of $c_2$, 
 we deduce that 
 $D_2$ lies in the interior of $\Cone {\Delta_1} B$; we can then choose an integral class $D$ in the interior of $\Cone {\Delta _1} {D_2}$.  Such a class $D$ satisfies $D^2 \cdot E_i >0$ for $i=1,2$, 
 and we therefore have a divisor $D$ with the required properties.

 We must finally deal with the case (in the above notation) where $E_1$ is known (with $E_1^3 \ge 0$)
but has $c_2 \cdot E_1 = 0$ (so $E_1$ and $-F_2$ generate the same ray), and therefore it is 
the unknown class $E_2$ 
which satisfies $c_2 \cdot E_2 >0$.  Note that our assumptions imply in this case that $E_1 ^3 >0$ since $E_1$ lying on the ray generated by $-F_2$ implies that it does not lie on the ray generated by $F_1$.  
We recall that $D_2^2$ is a positive multiple of $c_2$ and is therefore zero on $E_1$ and 
positive on $E_2$.  
If the ray generated by $D_2$ is \it rational\rm , we can choose $D_2$ to be integral and then 
$D = D_2$ is the divisor class we seek in this subcase.

The only remaining case is that of the previous paragraph but where 
 the ray generated by $D_2$ is not rational; here we shall need a further  idea.  We 
have reduced to the case where the class $E_1$ is known and has $c_2 \cdot E_1 =0$ and $E_1 ^3 >0$; 
 for any very ample divisor $H$ on $X$ generic in its linear system, 
we noted in Lemma 2.3  that the curve $C = E_1 |_H$ on $H$ then has the property that the pair $(H, \mu C)$ is klt for all $\mu < {5\over 6}$.
  Choose an integral $D_0$ in the interior of $\Cone {D_2} B$, and then  we may assume that 
a specific $D_2$ has been chosen in its given ray of the form  
$D_2 = D_0 + \lambda E_1$, with $\lambda$ a positive non-rational real number; for notational simplicity we now set $D = D_2$.  For $m$ a positive integer, we define the \it round-up \rm of $mD$ to be 
the integral divisor $\lceil {mD} \rceil = mD_0 + \lceil {m \lambda} \rceil E_1$.  We say that $m$ is \it 
relevant \rm if 
$\lceil m\lambda \rceil  - m\lambda <{5\over 6}$; even when  $\lambda$ is irrational this last condition holds for infinitely many $m$ (depending  only on $\lambda$).

For a given relevant integer $m$, we can choose a very ample linear system $|H|$ on $X$ for which 
$h^2 (X, \O _X ( \lceil {mD} \rceil +H )) =0$, and for 
$H$ generic in this linear system, Lemmas 3.3 and  4.3 imply that  
$D|_H$ is nef and big on $H$.  

We first claim that $h^1 (H, \O_H (K_H + \lceil {mD} \rceil )) =0$ for the given relevant $m$. 
 The claim follows from  the Kawamata--Viehweg version of Kodaira Vanishing for real divisors on the surface $H$, since $D|_H$ is big and nef on $H$ and $(H, \mu C)$ is klt for $\mu < {5\over 6}$.  Usually this form of Kodaira Vanishing 
is stated with $\Q$-divisors, but it may be extended to real divisors by the argument in \cite{Laz}, Vol. II, Remark 9.1.23, and the 
precise version used is stated in  \cite{Bir}, Theorem 7.4.  Given that $h^2 (X, \O _X ( \lceil {mD} \rceil +H )) =0$, we 
 may use  the short exact sequence 
$$ 0 \to \O_X (\lceil {mD} \rceil ) \to \O_X ( \lceil {mD} \rceil + H) \to \O_H (K_H + \lceil {mD} \rceil ) \to 0,$$ to deduce that $h^2 (X, \O _X ( \lceil {mD} \rceil )) =0$.  This is therefore true for all relevant $m$.

Now choose such an $m$ with  $(\lceil {mD} \rceil )^3 >0$  suitably large, i.e. so that 
$$\chi (X, \O_X (\lceil {mD} \rceil )) = {1\over 6} (\lceil {mD} \rceil )^3 + {1\over {12}} c_2 \cdot \lceil {mD} \rceil  >1,$$and hence $h^0 (X, \O_X (\lceil {mD} \rceil )) >1$ for some large (relevant) $m$ 
depending only on the known data.  Thus in this case, we take $M 
= \lceil mD \rceil$.

 (b)  We now consider case (2) when the cubic form is as in Case (b) from Section 5.  
 Thus (rational) coordinates may be chosen so that the cubic takes the form $x^2 y$.  Here positivity is just given by $y>0$ ($x\ne 0$) and the index condition is $x \ne 0$.  
Without loss of generality we may take $ P = \Cone A B$ with $A= (0,1)$ and $B = (1,0)$.  Furthermore we may assume that in 
fact $A$ and $B$ are integral and primitive, at the expense of the cubic form being a rational multiple of $x^2 y$.  In this case $A$ and $B$ will generate a sublattice  of 
$H^2 (X, \Z)$ (depending only on the cubic form) of finite index.   If $c_2$ were to vanish on a wall of $P$, then this wall would contain one of the rigid non-movable surfaces $E_i$, 
 and this is ruled out by Lemma 5.4.   Thus our assumptions imply 
 that $c_2$ is strictly positive on $P$, as clearly it cannot cut $P^\circ$.  As in 
 Proposition 6.1, we may assume $E_1^3 >0$ and $E_2 ^3 <0$
 with $E_1 \in \Cone {-B} A$ and $E_2 \in \Cone B {-A}$.
In particular the class of $E_1$ is determined up to finitely many possibilities by Proposition 2.2,
and so we may assume the class of $E_1$ is given.  

 As in part (a), we consider a class $D_2 \in P^\circ$, unique up to positive multiples, for which $c_2$ is a positive multiple of $D_2 ^2$.  Here there is a considerable simplification compared to part (a), in that $D_2$ may be taken to be rational.   To see this let $F_2$ be a rational generator of the line $c_2 =0$ in the quadrant $\Cone B {-A}$, which in coordinates may be written as  $F_2 = (p, -q)$ for some positive rationals $p,q$.  The quadratic equation $F_2 \cdot Z^2 =0$ is satisfied by $Z=A$, and so the other solution $D_2 \in P^\circ$ may also be taken to be rational --- in coordinates we can explicitly take $D_2 = (2p, q)$. We can then take $D_2$ to be 
integral and primitive (hence the class itself can be uniquely defined).   
Since $D_2^2$ is a positive multiple of $c_2$, we have 
$D_2 ^2 \cdot E_i \ge 0$ for $i=1, 2$ (at least one inequality being strict), and then  $D=D_2$ is the divisor class we seek in this subcase.

 (c)  Finally we consider case (2) when the cubic form is as in Case (c) from Section 5; we adopt the conventions as at the start 
 of Section 5, with $P= \Cone A B$ where $A = (1, 1/\sqrt{3})$ and $B= (1,0)$.
 We know that neither $E_i$ 
 lies in $P$ by Lemma 5.4.  We then have the possibility that $c_2=0$ is the line $y = -x/\sqrt 3$, i.e. gives a  
 rational solution to the Hessian quadratic.  By the Hodge index theorem, the index of the quadratic form 
 corresponding to each $E_i$ is either $(1,1)$ or $(1,0)$, and so the Hessian form at each $E_i$ must take non-positive values.
 Moreover our assumptions then imply that one of the $E_i$ (without loss of generality $E_1$)
 has $c_2 \cdot E_1 =0$.
 In this case, the wall $y = x/\sqrt 3$ of $P$ (generated by $A$)  is in fact generated by an integral divisor $\Delta$ with $\Delta \cdot E_1 \equiv 0$, some multiple of which defines the contraction of $E_1$ as in Proposition 5.5.  We observed in Remark 5.6 that we had boundedness of the threefolds in this case

 Therefore the only remaining case is when 
 the line $c_2 = 0$ divides the interior of the cone $\Cone {(1, 0)} {(1, - 1/\sqrt 3 )}$ and its negative, 
 since otherwise $c_2$ would be negative on one of the $E_i$.  We let $F_1 =
 (-1, 1/\sqrt 3 )$ and $F_2 $ in the interior of $ \Cone B {-F_1}$ generating the line $c_2 =0$, 
 and $D_2 \in P^\circ$ has been chosen with $D_2 ^2 \cdot F_2 =0$, i.e. $D_2^2$ is a positive multiple of $c_2$.  
 By consideration of the index, we may suppose then 
 that $E_1 \in \Cone {-F_2 }{ F_1}$ and $E_2 \in \Cone B {F_2}$, with $E_2$ not in the ray generated by $B$ from Lemma 5.4.  
 In particular $E_1 ^3 > 0$ and $E_2^3 <0$.  Thus there are only finitely many possibilities for the class of $E_1$; 
 for any given choice, we can choose 
  $\Delta_1 \in P^\circ$ with $\Delta_1 ^2 \cdot E_1 =0$, and we may choose an integral 
  class $D$ in the interior of $\Cone {\Delta _1} {D_2}$, for which therefore 
 $D^2 \cdot E_i >0$ for $i=1,2$.  Then $D$ is the integral divisor class we seek in this subcase.
 \end{proof}\end{thm}
 
  \section{Proof of Theorem 0.3 for case (3)}

Recall from  Section 5 that 
 case (3) was when there were two rigid non-movable surface classes, 
  $E_1$ and $E_2$ with the $c_2 ·\cdot E_i $ having different signs.  Throughout we can assume that $X$ is general in moduli.  Moreover we showed just after Proposition 5.1 that $\Effb (X) = \Cone {E_1} {E_2}$. 
  
  Because $c_2$ is negative on one of the $E_i$, using $c_2$ to  bound the 
 rigid non-movable fixed 
 part (as we did in Theorem 6.1 of the last section) does not  work, but we instead have to introduce a different method.  The crucial observation is that if say $c_2 \cdot E_1 <0$, the class 
 of $E_1$ is determined up to finitely many possibilities, and for each possibility there exists a real class $\Delta _1 \in P^\circ$ such that $\Delta _1 ^2 \cdot E_1 =0$.
 
 \medskip
 
\begin{thm}  The result in Theorem 0.3 is true in case (3) of Section 5.

\begin{proof}

Throughout the proof, we may assume that $X$ is general in moduli.  
We may also assume that the component $P^\circ$ of the positive index cone that contains the 
\Kahler cone has been specified, with closure denoted by $P$.
By assumption we have rigid non-movable surfaces $E_i$ with $i = 1,2$, where say $c_2 \cdot E_1 <0$ and $c_2 \cdot E_2 >0$.  Thus by Proposition 2.2, 
given the cubic and linear forms, there are only finitely many possible classes for $E_1$, and for each of these $E_1 ^3 >0$.  We may assume now that $E_1$ is known, and by Remark 5.7 that it is not a root of the Hessian quadratic.  Thus by Proposition 5.5, there is a class $\Delta _1 \in P^\circ$, unique up to positive multiples, such that $\Delta _1 ^2 \cdot E_1 =0$; it is the existence of such a $\Delta _1$ which makes the proof work.

(a) We consider first the case when the cubic has three distinct real 
linear factors, namely Case (a) from Section 5.  With our usual coordinates, we assume that $P =\Cone A B$, where $A= (0,1)$ and $B= (1,0)$. 
Since $E_1 ^3 >0$, we may assume without loss of generality 
(switching the two coordinates if necessary) that $E_1 \in \Cone {(-1,0)} {(-1, 1)}$.
 With $\Delta _1$ as defined above, we note that $\Delta _1 ^2 \cdot E_2 >0$ (since any ample divisor is a convex combination of $E_1$ and $E_2$, we would otherwise obtain a contradiction to 
  Lemma 3.3).  We deduce that $\Delta _1 |_H$ is nef on the   generic very ample divisor $H$ 
  by Lemma 4.3.  We note that $c_2 \cdot B >0$, since otherwise $c_2$ would be negative on $\calK \subset \Cone {E_1} B$.  We can therefore 
 choose integral $F$ in the interior of $ \Cone {E_1} B$ with $c_2 \cdot F =0$, i.e. $F$ is in the line $c_2 =0$.  Since the line $\Delta _1 ^2 =0$ contains $E_1$ by definition, and 
 $\Delta _1 ^2 \cdot Z >0$ for any $Z \in P^\circ$ by Lemma 3.3, we see that 
 $\Delta _1 ^2  \cdot B >0$ and hence $\Delta _1 ^2  \cdot F >0$.
     Moreover, since $c_2$ is non-negative on $\Movb (X)$
    by Remark 1.2, we have $\Movb (X) \subset \Cone F {-E_1}$.  However
    by Proposition 5.5, we can find an explicit integral $R$ in the interior of $ \Cone F {-E_1}$ such that $\Movb (X) \subset \Cone F R$; note that $\Delta _1 ^2 \cdot R >0$ (since 
    $\Delta _1 ^2 \cdot F >0$ and $\Delta _1 ^2 \cdot (-E_1) = 0$).
  
  Suppose first that $\Delta _1$ 
 may be chosen rational, and so may also be assumed integral; then 
 by Proposition 4.1, we can find a positive integer $m$ depending only on known data (including the choice of $E_1$ and $\Delta _1$)
 such that $|m \Delta _1 | = |L| + aE_1 + bE_2$, where some multiple of $L$ is mobile (cf. Remark 4.2) and $a,b$ non-negative integers.  
 Since $\Delta _1 ^2 \cdot E_2 >0$, we note 
 in particular that $\Delta _1 ^2 \cdot L \le m\Delta _1 ^3$, and so
 $\Delta _1 ^2 \cdot L$ is bounded above.   Since 
  $L \in \Cone {F} R$, with 
 $\Delta _1 ^2 \cdot F >0$ and $\Delta _1 ^2 \cdot R >0$ known values, and 
 $\Delta _1 ^2 \cdot L$ is explicitly bounded above, 
  there are only finitely many possibilities for the class of the movable divisor $L$ (since the lattice $\Z F + \Z R$ has finite index in $H^2 (X, \Z )$). Boundedness then follows from 
  Proposition 1.1 and Theorem 5.7 (i).

For the case when $\Delta _1$ cannot be chosen rational, we employ the argument from 
part (a) of the proof of Theorem 6.3.  For simplicity of notation, we write $D$ for $\Delta _1$.
For a generic very ample divisor $H$ on $X$, we again have  
the curve $C = E_1 |_H$ where by 
Lemma 2.3 the pair $(H, \mu C)$ is  klt for $\mu < {5\over 6}$; this fact is slightly easier here since we do not need to worry about the case $c_2 \cdot E_1 = 0$ 
and $E_1^3 >0$  in that Lemma.
With integral $D_0$  chosen in the interior of $\Cone D B$, we may assume that 
$D = D_0 + \lambda E_1$, with $\lambda$ a positive non-rational real number.  
The argument from part (a) of the proof of Theorem 6.3 then shows that $h^2 (X, \O _X ( \lceil {mD} \rceil )) =0$ for infinitely many (relevant) $m>0$.
  Now choose  $m$ large enough so that 
$$\chi (X, \O_X (\lceil {mD} \rceil )) = {1\over 6} (\lceil {mD} \rceil )^3 + {1\over {12}} c_2 \cdot \lceil {mD} \rceil  >1,$$ and hence we find $m$ depending 
only  the known data with $h^0 (X, \O_X (\lceil {mD} \rceil )) >1$.  Write 
$| \lceil {mD} \rceil  | = |L |+ a E_1 + bE_2$ where by Remark 4.2 some multiple of $L$ is mobile; note that  $D^2 \cdot \lceil {mD} \rceil = m D^3$ and so $D^2 \cdot L$ is also bounded above.  As 
$L \in \Cone {F} R$, the argument from the previous paragraph shows that 
there are only finitely many possibilities for the class of the movable  divisor $L$ and boundedness then follows from Proposition 1.1 and Theorem 5.7 (i). 

(b)  For the case when the cubic has a linear factor of multiplicity two, namely Case (b) from Section 5, the same arguments as for Case (a) work.
Since $c_2 \cdot E_1 <0$ implies that $E_1^3 >0$ by Proposition 2.2, we have that 
$E_1 , E_2$ are  now in the interiors of 
$\Cone {-B} A$ and $\Cone B {-A}$ respectively.
  Here $A = (0,1)$, $B= (1,0)$ and we have taken 
$P = \Cone A B$.
The integral class $E_1$ may be assumed known and we 
can then choose  $\Delta _1 \in P^\circ$ with $\Delta _1 ^2 \cdot E_1 =0$,
with the corresponding ray unique.  However, by the argument 
given in the second paragraph of part (b) of the proof of Theorem 6.3, in this case 
$\Delta _1$ may be chosen integral, since the quadratic $E_1 \cdot Z^2 =0$ is satisfied at the integral class $A$.  We note that $\Delta _1 ^2 \cdot E_2 >0$ (since any ample divisor is a convex combination of $E_1$ and $E_2$, we would otherwise obtain a contradiction to  Lemma 3.3).
Therefore we deduce that $\Delta _1 |_H$ is nef on the generic very ample divisor $H$.
Noting that $c_2 \cdot B >0$, since otherwise $c_2$ would be negative on $\calK \subset \Cone {E_1} B$,  we
can define the integral class $F$
in the interior of $ \Cone {E_1} B$ with $c_2 \cdot F =0$; since the line $\Delta _1 ^2 =0$ contains $E_1$ by definition, we observe 
  that $\Delta _1 ^2  \cdot F >0$. 
    Moreover, by Proposition 5.5, there exists an explicit integral class $R$ in the interior of $ \Cone F {-E_1}$ such that $\Movb (X) \subset \Cone F R$; note that $\Delta _1 ^2 \cdot R >0$ (since 
    $\Delta _1 ^2 \cdot F >0$ and $\Delta _1 ^2 \cdot (-E_1) = 0$).  
By Proposition 4.1, 
we may as in part (a) 
then find a positive integer $m$ depending only on known data such that $|m \Delta _1 | = |L| + aE_1 + bE_2$, where some multiple of $L$ is mobile (cf. Remark 4.2) 
and $a,b$ non-negative integers.  In particular we 
note  that $\Delta _1 ^2 \cdot L$ is bounded above.  
 Since $L \in \Cone {F} R$, with $\Delta _1 ^2 \cdot F >0$ and $\Delta _1 ^2 \cdot R >0$ 
 known values, 
 and $\Delta _1 ^2 \cdot L$ is explicitly bounded above, 
  there are only finitely many possibilities for the class of the movable divisor $L$. Boundedness then follows from Proposition 1.1 and Theorem 5.7 (i).

(c)  Finally we consider Case (c) for the cubic, as defined in Section 5.  With the usual coordinates,  
we have $P = \Cone A B$ where $A = (1, 1/\sqrt{3} )$ and $B= (1,0)$.  
Since $c_2 \cdot E_1 <0$ and thus $E_1^3 >0$ by Proposition 2.2, 
 index considerations then imply that $E_1 \in \Cone {-B} {(-1, 1/\sqrt 3 )}$; 
it is clearly not 
a multiple of $-B$ (since $B\cdot H^2 >0$ for all $H\in P^\circ$, or by using Proposition 5.1).    
As in part (c) of the 
proof of Theorem 6.3, we may also assume without loss of generality that $E_1$ is not in the ray 
generated by $(-1, 1/\sqrt{3})$ (this case is covered by Remark 5.6).  So $E_1$ may be assumed to be in the interior 
of the cone $\Cone {-B} {(-1, 1/\sqrt 3 )}$.
Since $-E_2 \not\in P$, index considerations imply that $E_2$ is in the interior of $\Cone {(1, 0)} {(1, -1/\sqrt 3 )}$; it is not in the ray generated by $(1, -1/\sqrt 3 )$, since then (as calculated in part (c) of the proof 
of Theorem 6.2)
we would have $E_2 \cdot H^2 < 0$ for any $H\in P^\circ$, and it is not in the ray generated by $B$ by Lemma 5.4.  For each of (the finite number of) possibilities for the class $E_1$, we can choose  $\Delta _1 \in P ^\circ$, unique up to positive multiples, 
 with $E_1 \cdot \Delta _1 ^2 = 0$.

  The ray generated by $\Delta _1 \in P^\circ$ may however be rational or irrational.
 We remark as in the previous two parts that $\Delta _1 ^2 \cdot E_2>0$ and then in the same way as there 
 we deduce that $\Delta _1 |_H$ is nef on a  generic very ample divisor $H$.  
 
 We choose integral classes $F$ and $R$ in exactly the same way as in part (a), where $\Delta _1 ^2 \cdot F >0$ and $\Delta _1 ^2 \cdot R >0$ and $\Movb (X) \subset \Cone F R$.
  When the ray generated by $\Delta _1$ is rational,  we may take $\Delta _1$  to be integral  
  and we can find an $m>0$ depending only on the known data with 
$h^0 (X , \O_X (m\Delta _1)) >1$; 
we may then conclude by the same argument as in  the second paragraph of part (a).   When the ray generated by $\Delta _1$ is not rational, the more delicate argument from part  (a) of this proof needs
  be employed;
in particular the last paragraph in the proof of part (a) translates unchanged to this case, to 
give a finite number of non-trivial integral classes $L$, at least one of which will have a multiple which is mobile.  Boundedness then follows from Proposition 1.1 and Theorem 5.7 (i).
\end{proof} \end{thm}

\smallskip

Thus Theorem 0.3 follows from Theorem 0.1, Theorem 6.2, Theorem 6.3 and Theorem 7.1.
\bibliographystyle{amsplain}

\end{document}